\RequirePackage{fix-cm}
\documentclass{article}

\usepackage{graphicx}%
\usepackage{multirow}%
\usepackage{amsmath,amssymb,amsfonts}%
\usepackage{amsthm}%
\usepackage{mathrsfs}%
\usepackage[title]{appendix}%
\usepackage{xcolor}%
\usepackage{textcomp}%
\usepackage{manyfoot}%
\usepackage{booktabs}%
\usepackage{algorithm}%
\usepackage{algorithmicx}%
\usepackage{algpseudocode}%
\usepackage{listings}%
\usepackage{xcolor}
\usepackage{url}
\usepackage{hyperref}
\hypersetup{
	colorlinks=true,
	linkcolor=blue,
	citecolor=blue,
	urlcolor=blue
}

\theoremstyle{thmstyleone}%

\theoremstyle{thmstyletwo}%
\newtheorem{remark}{Remark}%
\theoremstyle{thmstylethree}%
\title{Disjunctive programming for bi-objective combinatorial optimization problems with integer data and box uncertainty}
\author{Alejandro Crema\\
	Escuela de Computaci\'on, Facultad de Ciencias,\\
	Universidad Central de Venezuela, Caracas, Venezuela.\\
	email: alejandro.crema@ciens.ucv.ve,alejandro.crema@gmail.com}

\makeatletter
\renewcommand{\@biblabel}[1]{}
\makeatother

\begin{document}
	\maketitle
\abstract{In this paper we study bi-objective combinatorial optimization problems with integer data and box uncertainty. By employing extreme scenario characterizations on the uncertainty box, we define preorders on the feasible set whose induced posets allow us to identify coverages by finding minimal elements. We present exact and approximate parametric algorithms based on set-covering formulations to characterize the trade-off between cardinality and quality for coverages. Computational experience on bi-objective set covering, selection and shortest path problems is presented.}\\


\maketitle

\section{Introduction}

Let $n \geq 1$ and 
\begin{itemize}
	\item  let $\mathcal{X} \subseteq {\lbrace 0,1\rbrace}^n$ be a non-empty set that can be represented by a 0-1 mixed integer linear programming (0-1-MILP) formulation,
	\item let $\underline{\mathbf{a}},\overline{\mathbf{a}} \in \mathbb{Z}^n$ with $ \mathbf{0} \leq \underline{\mathbf{a}} \leq \overline{\mathbf{a}}$, let $[\underline{\mathbf{a}},\overline{\mathbf{a}}] = \lbrace \boldsymbol{a} \in \mathbb{Z}^n:\underline{\mathbf{a}} \leq \boldsymbol{a} \leq \overline{\mathbf{a}} \rbrace$,
	\item let $\underline{\mathbf{c}},\overline{\mathbf{c}} \in \mathbb{Z}^n$ 
	with $ \mathbf{0} \leq \underline{\mathbf{c}} \leq \overline{\mathbf{c}}$, 
	let $[\underline{\mathbf{c}},\overline{\mathbf{c}}] = \lbrace \boldsymbol{c} \in \mathbb{Z}^n:\underline{\mathbf{c}} \leq \boldsymbol{c} \leq \overline{\mathbf{c}} \rbrace$ and
	\item let $\mathcal{U} = [\underline{\mathbf{a}},\overline{\mathbf{a}}] \times [\underline{\mathbf{c}},\overline{\mathbf{c}}]$ be the {\it box uncertainty} set.
\end{itemize}

Let $\mathcal{P}(\boldsymbol{a},\boldsymbol{c})$ be a bi-objective combinatorial optimization (BOCO) problem in $\boldsymbol{x}$, parameterized  by $(\boldsymbol{a},\boldsymbol{c}) \in \mathcal{U}$, defined as follows:
\begin{flalign*}
	\min&\;\;\boldsymbol{a}^T\boldsymbol{x}, \boldsymbol{c}^T\boldsymbol{x} &\;&&\mathcal{P}(\boldsymbol{a},\boldsymbol{c})&\\
	s.t.&\;\boldsymbol{x} \in \mathcal{X}&\
\end{flalign*}

The $\mathcal{U}$-$BOCO$ problem is the set of $BOCO$ problems defined as:
\begin{flalign*}
	&\;\;\lbrace \mathcal{P}(\boldsymbol{a},\boldsymbol{c}): (\boldsymbol{a},\boldsymbol{c}) \in \mathcal{U} \rbrace &\;&&\mathcal{U}\mbox{-}BOCO &\
\end{flalign*}

Although the results and algorithms to be presented can easily be reformulated for any number of objectives, we have chosen to frame everything in terms of two objectives in order to simplify the notation and because all the experiments we describe were conducted using two objectives.\\

We say that  $\mathcal{C}(\boldsymbol{a},\boldsymbol{c})$ is an $(\boldsymbol{a},\boldsymbol{c})$-{\it coverage} if and only if (i) $\mathcal{C}(\boldsymbol{a},\boldsymbol{c})\subseteq \mathcal{X}$ and (ii) for all  $\boldsymbol{x} \in \mathcal{X}$ there exists $\boldsymbol{y} \in   \mathcal{C}(\boldsymbol{a},\boldsymbol{c})$ such that  $(\boldsymbol{a}^T\boldsymbol{y}, \boldsymbol{c}^T\boldsymbol{y}) \leq (\boldsymbol{a}^T\boldsymbol{x}, \boldsymbol{c}^T\boldsymbol{x})$. If  $\mathcal{C}(\boldsymbol{a},\boldsymbol{c})$ is any minimum-cardinality   $(\boldsymbol{a},\boldsymbol{c})$-coverage then  $\mathcal{PF}(\boldsymbol{a},\boldsymbol{c}) = \lbrace (\boldsymbol{z},\boldsymbol{w}): \exists \boldsymbol{y} \in \mathcal{C}(\boldsymbol{a},\boldsymbol{c}) \mbox{ with } (\boldsymbol{z},\boldsymbol{w}) = (\boldsymbol{a}^T\boldsymbol{y}, \boldsymbol{c}^T\boldsymbol{y}) \rbrace$ is the {\it Pareto frontier} for $(\boldsymbol{a},\boldsymbol{c})$. \\

For information on the importance of finding $\mathcal{PF}(\boldsymbol{a},\boldsymbol{c})$, as well as on algorithms for identifying it, one may consult  and \cite{Ehrgott2005} and the excellent reviews by \cite{Ehrgott2026}, \cite{halffmann2022}, and \cite{Lokman2026} along with their respective references.\\

An obvious generalization of the concept of coverage when considering uncertainty is the following: $\mathcal{C}$ is a $\mathcal{U}$-{\it coverage} if and only if (i) $\mathcal{C}\subseteq \mathcal{X}$ and (ii) for all $((\boldsymbol{a},\boldsymbol{c}),\boldsymbol{x}) \in \mathcal{U} \times \mathcal{X}$ there exists $\boldsymbol{y} \in \mathcal{C}$ such that $(\boldsymbol{a}^T\boldsymbol{y}, \boldsymbol{c}^T\boldsymbol{y}) \leq (\boldsymbol{a}^T\boldsymbol{x}, \boldsymbol{c}^T\boldsymbol{x})$. Let $(\boldsymbol{a},\boldsymbol{c}) \in \mathcal{U}$ and let $\mathcal{C}$ be a minimum-cardinality $\mathcal{U}$-coverage, then we will see that it is very easy to find a minimum-cardinality $(\boldsymbol{a},\boldsymbol{c})$-coverage  contained in $\mathcal{C}$. We present algorithms for finding a minimum-cardinality $\mathcal{U}$-coverage.\\

The cardinality of a minimum-cardinality $(\boldsymbol{a},\boldsymbol{c})$-coverage can be problematic; therefore, various approaches have been developed to identify a suitable representation. For information on the criteria used to determine whether a set constitutes a suitable representation, we refer to the works of \cite{Faulkenberg2010}, \cite{Humbert2026}, and \cite{Sayin2000}, as well as the references cited therein.\\

In this article, we apply an extension of one of the criteria for qualifying a set as a good representation to the case involving uncertainty. Let $(\boldsymbol{a},\boldsymbol{c}) \in \mathcal{U}$ and let $\epsilon \in \mathbb{N} \cup \lbrace 0 \rbrace$. $\mathcal{C}(\epsilon,\boldsymbol{a},\boldsymbol{c})$ is an $\epsilon$-$(\boldsymbol{a},\boldsymbol{c})$-{\it coverage} if and only if (i) $\mathcal{C}(\epsilon,\boldsymbol{a},\boldsymbol{c})\subseteq \mathcal{X}$ and (ii) for all $\boldsymbol{x} \in \mathcal{X}$ there exists $\boldsymbol{y} \in \mathcal{C}(\epsilon,\boldsymbol{a},\boldsymbol{c})$ such that $(\boldsymbol{a}^T\boldsymbol{y}, \boldsymbol{c}^T\boldsymbol{y}) \leq (\boldsymbol{a}^T\boldsymbol{x}+\epsilon, \boldsymbol{c}^T\boldsymbol{x}+\epsilon)$. An obvious generalization of the concept of $\epsilon$-coverage when considering uncertainty is the following:  $\mathcal{C}(\epsilon)$ is an  $\epsilon$-$\mathcal{U}$-{\it coverage} if and only if (i)  $\mathcal{C}(\epsilon)\subseteq \mathcal{X}$ and  (ii) for all  $((\boldsymbol{a},\boldsymbol{c}),\boldsymbol{x}) \in \mathcal{U} \times \mathcal{X}$ there exists $\boldsymbol{y} \in \mathcal{C}(\epsilon)$ such that $(\boldsymbol{a}^T\boldsymbol{y}, \boldsymbol{c}^T\boldsymbol{y}) \leq (\boldsymbol{a}^T\boldsymbol{x}+\epsilon, \boldsymbol{c}^T\boldsymbol{x}+\epsilon)$.\\

We present a {\it parametric} algorithm to find, starting from a minimum-cardinality $\mathcal{U}$-coverage, a minimum-cardinality $\epsilon$-$\mathcal{U}$-coverage for all $\epsilon \in \mathbb{N} \cup \lbrace 0 \rbrace$. Let $\epsilon \in \mathbb{N} \cup \lbrace 0 \rbrace$; we present an algorithm to find an $\epsilon$-$\mathcal{U}$-coverage without the use of  a minimum-cardinality $\mathcal{U}$-coverage; in this case, it is possible to use an approximate parameterization to find a $(\delta+\epsilon)$-$\mathcal{U}$-coverage for all $\delta \in \mathbb{N} \cup \lbrace 0 \rbrace$.\\

The underlying idea in this work is that, once the {\it true} scenario, denoted as $(\boldsymbol{a},\boldsymbol{c})$, is revealed, it is possible to easily obtain a minimum-cardinality $(\boldsymbol{a},\boldsymbol{c})$-coverage, a minimum-cardinality $\epsilon$-$(\boldsymbol{a},\boldsymbol{c})$-coverage or simply an $\epsilon$-$(\boldsymbol{a},\boldsymbol{c})$-coverage. This concept is analogous to the one employed in the work of \cite{Buchheim2017}, as well as in numerous subsequent studies that cite it, in which algorithms are presented that compute $k$ solutions with the goal of selecting the best one once the true scenario is revealed. Even if that idea is not valid for the decision-maker, knowing a minimum-cardinality $\mathcal{U}$-coverage, a minimum-cardinality $\epsilon$-$\mathcal{U}$-coverage, or simply an $\epsilon$-$\mathcal{U}$-coverage can still be useful for addressing uncertainty through other approaches.\\

Some of the fundamental concepts of robust programming have been extended to the multi-objective case. Key ideas can be found in the works of \cite{Ehrgott2014}, \cite{Kuhn2016}, and \cite{Ide2016}. One of the proposed approaches is as follows: the problem under uncertainty is reduced to a deterministic, non-linear multiple objective problem, in which each objective is replaced by its worst-case evaluation, specifically the problem considered, with two objectives, is defined as follows:
\begin{flalign*}
	\min&\;\;\max_{\boldsymbol{a} \in \mathcal{U}_a} \boldsymbol{a}^T \boldsymbol{x},\max_{\boldsymbol{c} \in \mathcal{U}_c} \boldsymbol{c}^T \boldsymbol{x}&\;&&\widehat{\mathcal{P}}&\\
	s.t.&\;\boldsymbol{x} \in \mathcal{X}&\
\end{flalign*}
 
with $\mathcal{U}_a$ and $\mathcal{U}_c$ the uncertainty sets for $\boldsymbol{a}$ and  $\boldsymbol{c}$. In the case of box-type uncertainty, this approach leads to obtaining $\mathcal{PF}(\overline{\mathbf{a}},\overline{\mathbf{c}})$. The work of \cite{Raith2018} considers uncertainty sets defined by cardinality constraints. Other recent approaches for addressing multi-objective problems under uncertainty can be found in \cite{Kentache2026}, with an arbitrary finite uncertainty set,  and \cite{Gutierrez2026}, with countably many scenarios.\\

In Section 2 we present basic algebraic results and results from the theory of relations on finite sets that allow us to find a minimum-cardinality $\mathcal{U}$-coverage by finding the set of minimal elements of an appropriate partially ordered set defined on $\mathcal{X}$ according to $\mathcal{U}$.\\

In Section 3 we present the algorithms. In 3.1 we present two algorithms to find, starting from $\boldsymbol{x} \in \mathcal{X}$, a solution $\boldsymbol{y} \in \mathcal{X}$ that can be included in a minimum-cardinality $\mathcal{U}$-coverage. These algorithms are used to define algorithms to be presented in other sections.\\

In 3.2 we present a disjunctive programming algorithm  to find, given $\epsilon \in \mathbb{N} \cup \lbrace 0 \rbrace$, an $\epsilon$-$\mathcal{U}$-coverage. If $\epsilon=0$ the algorithm finds a minimum-cardinality $\mathcal{U}$-coverage. In each iteration the algorithm finds a solution to be included in the coverage and adds a new binary variable and two constraints to be used in the next iterations. Our algorithms for the deterministic case, which can be found in \cite{Sylva2004} and \cite{Sylva2007}, were presented with these same characteristics.\\

In the deterministic case, various approaches have been developed to avoid the introduction of new binary variables; this generally allows for the identification of the Pareto frontier in problems with dimensions far exceeding those that can be effectively addressed using standard disjunctive programming. For information on these algorithms, see the works of \cite{Bektas2018}, \cite{Fleuren2024}, \cite{Sayin2014}, and \cite{mavrotas2013}. Drawing inspiration from the concepts presented in the work of \cite{Bektas2018}, we have developed an algorithm  to find an $\epsilon$-$\mathcal{U}$-coverage, presented in Section 3.3, that does not incorporate any additional variables. If $\epsilon=0$ the algorithm finds a minimum-cardinality $\mathcal{U}$-coverage.\\ 

The algorithms presented in 3.2 and 3.3 have the property that their output is only meaningful if they complete their tasks in their entirety. In 3.4, we present an {\it anytime} version that provides an $\epsilon$-$\mathcal{U}$-coverage if it is terminated at any given iteration, whether due to a time limit or some other stopping criterion, where $\epsilon$ is part of the solution. The algorithm presented in 3.4 is based on standard disjunctive programming and incorporates one binary variable and two constraints for each generated solution. Some remarks about the disjunctive problems are presented in Section 3.5.\\

In Section 3.6, we present exact and approximate parametric algorithms based on set-covering problems associated with a minimum-cardina-\
lity $\mathcal{U}$-coverage (or an $\epsilon$-$\mathcal{U}$-coverage). Section 4 presents computational results for $\mathcal{U}$-$BOCO$ problems, specifically the bi-objective set covering problem, the bi-objective selection problem, and the bi-objective shortest path problem. We will see how the algorithms presented in Sections 3.2 and 3.3 prove to be the best on different instances, with neither one dominating the other. Furthermore, we provide evidence from the parametric analysis that the anytime algorithm should not be dismissed, even though it is clearly inferior. Finally, conclusions are presented in Section 5.\\

Overall, combining  order structures with disjunctive programming and set-covering parameterizations provides a flexible framework for managing objective uncertainty in $\mathcal{U}$-$BOCO$ problems, while allowing exact or approximate coverage representations to be readily exploited once the true scenario is revealed.\\

A few words about the notation. If $\bullet$ is an optimization problem then $F(\bullet)$ is its feasible solution set and $v(\bullet)$ is its optimal value (if it exists). The vector with all components equal to 0(1) will be denoted $\mathbf{0}(\mathbf{1})$ and the dimensions must be clear from the context. For any $p \in \mathbb{N}$, $[p]$ denotes the index set $\lbrace 1,\cdots,p\rbrace$.\\

\section{Theoretical Results}
Below, we present basic algebraic results and results related to the theory of relations on finite sets, which will be used to define the algorithms.\\

 \subsection{Basic results}

Let $\underline{\mathbf{d}},\overline{\mathbf{d}} \in \mathbb{Z}^n$ with $ \mathbf{0} \leq \underline{\mathbf{d}} \leq \overline{\mathbf{d}}$ and let $[\underline{\mathbf{d}},\overline{\mathbf{d}}] = \lbrace \boldsymbol{d} \in \mathbb{Z}^n:\underline{\mathbf{d}} \leq \boldsymbol{d} \leq \overline{\mathbf{d}} \rbrace$. Let $\boldsymbol{d}^+(\boldsymbol{x})_j = \underline{\mathbf{d}}_j \boldsymbol{x}_j + \overline{\mathbf{d}}_j ( 1-\boldsymbol{x}_j)$ for all $j \in [n]$ and for all $\boldsymbol{x} \in \mathcal{X}$. Let $\boldsymbol{x} \in \mathcal{X}$; then $\boldsymbol{d}^+(\boldsymbol{x})$ is known as its \emph{most favorable scenario}(\cite{Aissi2009}).\\

In \textbf{Remarks} 1 through 4, we present elementary results related to $\boldsymbol{d}^+(\boldsymbol{x})$. Some proofs are omitted because the results follow immediately from the definition. Proofs are included for the non-trivial and less well-known results.\\

\begin{remark} \label{remark:remarkbasica}
	$\;$ 
	\begin{enumerate}
		\item $\boldsymbol{d}^+(\boldsymbol{x}) \in [\underline{\mathbf{d}},\overline{\mathbf{d}}]$ for all $\boldsymbol{x} \in \mathcal{X}$,
		\item $\boldsymbol{d}^+(\boldsymbol{x})^T \boldsymbol{x}=\underline{\mathbf{d}}^T \boldsymbol{x}$ for all $\boldsymbol{x} \in \mathcal{X}$,
		\item $\boldsymbol{d}^+(\boldsymbol{x})^T (\boldsymbol{u} - \boldsymbol{x}) \geq \boldsymbol{d}^T(\boldsymbol{u}- \boldsymbol{x})$ for all $\boldsymbol{x},\boldsymbol{u} \in \mathcal{X}$ and for all $\boldsymbol{d} \in [\underline{\mathbf{d}},\overline{\mathbf{d}}]$,
		\item $\boldsymbol{d}^+(\boldsymbol{y})^T\boldsymbol{x} \leq \underline{\mathbf{d}}^T \boldsymbol{y}$ if and only if $\boldsymbol{d}^T \boldsymbol{x} \leq \boldsymbol{d}^T \boldsymbol{y}$ for all $\boldsymbol{d} \in [\underline{\mathbf{d}},\overline{\mathbf{d}}]$ and for all $\boldsymbol{x},\boldsymbol{y} \in \mathcal{X}$.
	\end{enumerate}
\end{remark}

\textbf{Proof}:
\begin{itemize}
	\item Parts 1, 2, and 3 follow directly from the definition of $\boldsymbol{d}^+(\boldsymbol{x})$.
	\item For part 4:
	\begin{itemize}
		\item Let us suppose that $\boldsymbol{d}^+(\boldsymbol{y})^T\boldsymbol{x} \leq \underline{\mathbf{d}}^T \boldsymbol{y} = \boldsymbol{d}^+(\boldsymbol{y})^T\boldsymbol{y}$. We have that
		\[\mathbf{0} \geq \boldsymbol{d}^+(\boldsymbol{y})^T(\boldsymbol{x}-\boldsymbol{y}) \underset{\square}{\geq} \boldsymbol{d}^T(\boldsymbol{x}- \boldsymbol{y})\]
		for all $\boldsymbol{d} \in [\underline{\mathbf{d}},\overline{\mathbf{d}}]$, with $\square$ following from 3.
		\item Let us suppose that $\boldsymbol{d}^T \boldsymbol{x} \leq \boldsymbol{d}^T \boldsymbol{y}$ for all $\boldsymbol{d} \in [\underline{\mathbf{d}},\overline{\mathbf{d}}]$ and let $\boldsymbol{d} = \boldsymbol{d}^+(\boldsymbol{y})$, then $\boldsymbol{d}^+(\boldsymbol{y})^T\boldsymbol{x} \leq \boldsymbol{d}^+(\boldsymbol{y})^T\boldsymbol{y} =\underline{\mathbf{d}}^T \boldsymbol{y}\;\blacksquare$
	\end{itemize}
\end{itemize}

\begin{remark}\label{remark:remarkdxequaldy} 
	Let $\boldsymbol{x},\boldsymbol{y} \in \mathcal{X}$ then: $\boldsymbol{d}^+(\boldsymbol{y})^T\boldsymbol{x} \leq \underline{\mathbf{d}}^T \boldsymbol{y}$ and $\boldsymbol{d}^+(\boldsymbol{x})^T\boldsymbol{y} \leq \underline{\mathbf{d}}^T \boldsymbol{x}$ if and only if $\boldsymbol{d}^T\boldsymbol{x} = \boldsymbol{d}^T\boldsymbol{y}$ for all $\boldsymbol{{d}} \in [\underline{\mathbf{d}},\overline{\mathbf{d}}]$.
\end{remark}

\textbf{Proof}:
\begin{itemize}
	\item Let us suppose that $\boldsymbol{d}^T\boldsymbol{x} = \boldsymbol{d}^T\boldsymbol{y}$ for all $\boldsymbol{{d}} \in [\underline{\mathbf{d}},\overline{\mathbf{d}}]$, then we have:
	$\boldsymbol{d}^+(\boldsymbol{y})^T\boldsymbol{x} = \boldsymbol{d}^+(\boldsymbol{y})^T \boldsymbol{y} = \underline{\mathbf{d}}^T \boldsymbol{y}$ and $\boldsymbol{d}^+(\boldsymbol{x})^T \boldsymbol{y} = \boldsymbol{d}^+(\boldsymbol{x})^T \boldsymbol{x} = \underline{\mathbf{d}}^T \boldsymbol{x}$.
	\item Let us suppose that $\boldsymbol{d}^+(\boldsymbol{y})^T\boldsymbol{x} \leq \underline{\mathbf{d}}^T \boldsymbol{y}$ and $\boldsymbol{d}^+(\boldsymbol{x})^T\boldsymbol{y} \leq \underline{\mathbf{d}}^T \boldsymbol{x}$, then:
	\begin{enumerate}
		\item $\mathbf{0} \geq \boldsymbol{d}^+(\boldsymbol{y})^T \boldsymbol{x} - \underline{\mathbf{d}}^T \boldsymbol{y} = \boldsymbol{d}^+(\boldsymbol{y})^T \boldsymbol{x} - \boldsymbol{d}^+(\boldsymbol{y})^T \boldsymbol{y} \underset{\square}{\geq} \boldsymbol{d}^T\boldsymbol{x} - \boldsymbol{d}^T\boldsymbol{y}\;\;\forall \boldsymbol{d} \in [\underline{\mathbf{d}},\overline{\mathbf{d}}] $
		\item $\mathbf{0} \geq \boldsymbol{d}^+(\boldsymbol{x})^T \boldsymbol{y} - \underline{\mathbf{d}}^T \boldsymbol{x} = \boldsymbol{d}^+(\boldsymbol{x})^T \boldsymbol{y} - \boldsymbol{d}^+(\boldsymbol{x})^T \boldsymbol{x} \underset{\square}{\geq} \boldsymbol{d}^T\boldsymbol{y} - \boldsymbol{d}^T\boldsymbol{x}\;\;\forall \boldsymbol{d} \in [\underline{\mathbf{d}},\overline{\mathbf{d}}] $
	\end{enumerate}	
	with $\square$ from \textbf{Remark \ref{remark:remarkbasica}}. From (1) and (2) we have that $\boldsymbol{d}^T\boldsymbol{x} = \boldsymbol{d}^T\boldsymbol{y}$ for all $\boldsymbol{d} \in [\underline{\mathbf{d}},\overline{\mathbf{d}}]\;\blacksquare$
\end{itemize}

\begin{remark}\label{remark:remarkJxy}
	Let $J = \lbrace j \in [n]: \underline{\mathbf{d}}_j < \overline{\mathbf{d}}_j \rbrace$, let $\boldsymbol{x},\boldsymbol{y} \in \mathcal{X}$ and let us suppose that $J \neq \emptyset$ then:
	\begin{enumerate}
		\item If $\boldsymbol{d}^+(\boldsymbol{y})^T\boldsymbol{x} \leq \underline{\mathbf{d}}^T \boldsymbol{y}$ and $\boldsymbol{d}^+(\boldsymbol{x})^T\boldsymbol{y} \leq \underline{\mathbf{d}}^T \boldsymbol{x}$ then $\boldsymbol{x}_j = \boldsymbol{y}_j$ for all $j \in J$.
		\item If $\boldsymbol{d}^+(\boldsymbol{y})^T\boldsymbol{x} \leq \underline{\mathbf{d}}^T \boldsymbol{y}$, $\boldsymbol{d}^+(\boldsymbol{x})^T\boldsymbol{y} \leq \underline{\mathbf{d}}^T \boldsymbol{x}$ and $J = [n]$, then $\boldsymbol{x} = \boldsymbol{y}$.
	\end{enumerate}
\end{remark}

\textbf{Proof}:
\begin{enumerate}
	\item Let us suppose that $\boldsymbol{d}^+(\boldsymbol{y})^T\boldsymbol{x} \leq \underline{\mathbf{d}}^T \boldsymbol{y}$ and $\boldsymbol{d}^+(\boldsymbol{x})^T\boldsymbol{y} \leq \underline{\mathbf{d}}^T \boldsymbol{x}$. Let $j \in J$, let $\boldsymbol{d}^1,\boldsymbol{d}^2 \in [\underline{\mathbf{d}},\overline{\mathbf{d}}]$ such that $\boldsymbol{d}^1_k = \boldsymbol{d}^2_k$ for all $k \in [n]-\lbrace j \rbrace$, let $\boldsymbol{d}^1_j = \underline{\mathbf{d}}_j$ and let $\boldsymbol{d}^2_j = \overline{\mathbf{d}}_j$. From \textbf{Remark \ref{remark:remarkdxequaldy}} we have that ${\boldsymbol{d}^1}^T (\boldsymbol{x}-\boldsymbol{y})=0$ and ${\boldsymbol{d}^2}^T (\boldsymbol{x}-\boldsymbol{y})=0$, therefore $(\overline{\mathbf{d}}_j-\underline{\mathbf{d}}_j)(\boldsymbol{x}_j-\boldsymbol{y}_j)=0$ and then $\boldsymbol{x}_j = \boldsymbol{y}_j$.
	\item Follows directly from 1 $\blacksquare$
\end{enumerate}

\begin{remark}\label{remark:bounds}
	Let $\boldsymbol{\hat{x}},\boldsymbol{y} \in \mathcal{X}$ then
	\[\underline{\mathbf{d}}^T \mathbf{1} \underset{1}{\geq} \underset{\boldsymbol{x} \in \mathcal{X}}{\max}\; \underline{\mathbf{d}}^T \boldsymbol{x} - \underline{\mathbf{d}}^T \boldsymbol{y} \underset{2}{\geq} \underset{\boldsymbol{x} \in \mathcal{X}}{\max}\; \underline{\mathbf{d}}^T \boldsymbol{x} - \underset{\boldsymbol{x} \in \mathcal{X}}{\min}\;\boldsymbol{d}^+(\boldsymbol{x})^T \boldsymbol{y} \underset{3}{\geq}\]
	\[\underset{\boldsymbol{x} \in \mathcal{X}}{\max}\;( \underline{\mathbf{d}}^T \boldsymbol{x}- \boldsymbol{d}^+(\boldsymbol{x})^T \boldsymbol{y}) \underset{4}{\geq}
	 \underline{\mathbf{d}}^T \boldsymbol{\hat{x}}- \boldsymbol{d}^+(\boldsymbol{\hat{x}})^T \boldsymbol{y}\]
\end{remark}

\textbf{Proof}:
\begin{enumerate}
	\item Since $\boldsymbol{x},\boldsymbol{y} \in {\lbrace 0,1 \rbrace}^n$ and $\underline{\mathbf{d}} \geq \mathbf{0}$.
	\item Since $\boldsymbol{d}^+(\boldsymbol{x}) \ge \underline{\mathbf{d}}$ for all $\boldsymbol{x} \in \mathcal{X}$.
	\item Let $\boldsymbol{u} \in \mathcal{X}$ such that $\underset{\boldsymbol{x} \in \mathcal{X}}{\max}\;( \underline{\mathbf{d}}^T \boldsymbol{x}- \boldsymbol{d}^+(\boldsymbol{x})^T \boldsymbol{y}) = \underline{\mathbf{d}}^T \boldsymbol{u}- \boldsymbol{d}^+(\boldsymbol{u})^T \boldsymbol{y}$. Since $\boldsymbol{u} \in \mathcal{X}$, we have that $\underline{\mathbf{d}}^T \boldsymbol{u}- \boldsymbol{d}^+(\boldsymbol{u})^T \boldsymbol{y} \leq \underset{\boldsymbol{x} \in \mathcal{X}}{\max}\; \underline{\mathbf{d}}^T \boldsymbol{x} - \underset{\boldsymbol{x} \in \mathcal{X}}{\min}\;\boldsymbol{d}^+(\boldsymbol{x})^T \boldsymbol{y}$.
	\item Since $\boldsymbol{\hat{x}} \in \mathcal{X}$. $\blacksquare$
\end{enumerate}

\subsection{Relations on $\mathcal{X}$}
For standard definitions and basic properties regarding preordered sets, induced equivalence relations, quotient sets, and minimal elements in finite partially ordered sets, we refer the reader to \cite{Davey2002} and
\cite{Birkhoff1967}.\\

Let $\boldsymbol{a}^+(\boldsymbol{x})_j = \underline{\mathbf{a}}_j \boldsymbol{x}_j + \overline{\mathbf{a}}_j ( 1-\boldsymbol{x}_j)$ for all $j \in [n]$ and for all $\boldsymbol{x} \in \mathcal{X}$, and let $\boldsymbol{c}^+(\boldsymbol{x})_j = \underline{\mathbf{c}}_j \boldsymbol{x}_j + \overline{\mathbf{c}}_j ( 1-\boldsymbol{x}_j)$ for all $j \in [n]$ and for all $\boldsymbol{x} \in \mathcal{X}$.\\

Let $\epsilon \in \mathbb{N} \cup \lbrace 0 \rbrace$. Next we define the relations $\prec_+$, $\prec^+$, $\prec^\epsilon_+$ and $\prec^+_\epsilon$ on $\mathcal{X}$. Let $\boldsymbol{x},\boldsymbol{y} \in \mathcal{X}$ then:\\

$\boldsymbol{x} \prec_+ \boldsymbol{y}$ if and only if $\boldsymbol{a}^+(\boldsymbol{y})^T\boldsymbol{x} \leq \underline{\mathbf{a}}^T \boldsymbol{y}$,

$\boldsymbol{x} \prec^+ \boldsymbol{y}$ if and only if $\boldsymbol{c}^+(\boldsymbol{y})^T\boldsymbol{x} \leq \underline{\mathbf{c}}^T \boldsymbol{y}$,

$\boldsymbol{x} \prec^\epsilon_+ \boldsymbol{y}$ if and only if $\boldsymbol{a}^+(\boldsymbol{y})^T\boldsymbol{x} \leq \underline{\mathbf{a}}^T \boldsymbol{y}+\epsilon$  and

$\boldsymbol{x} \prec^+_\epsilon \boldsymbol{y}$ if and only if $\boldsymbol{c}^+(\boldsymbol{y})^T\boldsymbol{x} \leq \underline{\mathbf{c}}^T \boldsymbol{y}+\epsilon$.\\

Note that: 

$\boldsymbol{x} \nprec_+ \boldsymbol{y}$ if and only if $\boldsymbol{a}^+(\boldsymbol{y})^T\boldsymbol{x} \geq \underline{\mathbf{a}}^T \boldsymbol{y}+1$,

$\boldsymbol{x} \nprec^+ \boldsymbol{y}$ if and only if $\boldsymbol{c}^+(\boldsymbol{y})^T\boldsymbol{x} \geq \underline{\mathbf{c}}^T \boldsymbol{y}+1$,

$\boldsymbol{x} \nprec^\epsilon_+ \boldsymbol{y}$ if and only if $\boldsymbol{a}^+(\boldsymbol{y})^T\boldsymbol{x} \geq \underline{\mathbf{a}}^T \boldsymbol{y}+\epsilon+1$ and

$\boldsymbol{x} \nprec^+_\epsilon \boldsymbol{y}$ if and only if $\boldsymbol{c}^+(\boldsymbol{y})^T\boldsymbol{x} \geq \underline{\mathbf{c}}^T \boldsymbol{y}+\epsilon+1$.\\

From the definitions we have that $\prec_+$, $\prec^+$, $\prec^\epsilon_+$ and $\prec^+_\epsilon$ are reflexive relations, and from \textbf{Remark} \ref{remark:remarkbasica} we have that $\prec_+$ and $\prec^+$ are transitive relations.\\

Let $\prec$ be defined on $\mathcal{X}$ as follows: $\boldsymbol{x} \prec \boldsymbol{y}$ if and only if $\boldsymbol{x} \prec_+ \boldsymbol{y}$ and $\boldsymbol{x} \prec^+ \boldsymbol{y}$ for all $\boldsymbol{x},\boldsymbol{y} \in \mathcal{X}$. Let $\prec_\epsilon$ be defined on $\mathcal{X}$ as follows: $\boldsymbol{x} \prec_\epsilon \boldsymbol{y}$ if and only if $\boldsymbol{x} \prec^\epsilon_+ \boldsymbol{y}$ and $\boldsymbol{x} \prec^+_\epsilon \boldsymbol{y}$ for all $\boldsymbol{x},\boldsymbol{y} \in \mathcal{X}$.\\

We have that $\prec$ is a reflexive and transitive relation, that is, $(\mathcal{X},\prec)$ is a preordered set. Let $\sim$ be a relation defined on $\mathcal{X}$ as follows: $\boldsymbol{x} \sim \boldsymbol{y}$ if and only if $\boldsymbol{x} \prec \boldsymbol{y}$ and $\boldsymbol{y} \prec \boldsymbol{x}$ for all $\boldsymbol{x},\boldsymbol{y} \in \mathcal{X}$. We have that $\sim$ is an equivalence relation on $\mathcal{X}$. If $\boldsymbol{x} \in \mathcal{X}$, let $[\boldsymbol{x}]$ be its equivalence class and let $\frac{\mathcal{X}}{\sim} = \lbrace [\boldsymbol{x}]:\boldsymbol{x} \in \mathcal{X} \rbrace$, that is, $\frac{\mathcal{X}}{\sim}$ is the quotient set.\\

Let $\preceq$ be a relation on $\frac{\mathcal{X}}{\sim}$ defined as follows: $[\boldsymbol{x}] \preceq [\boldsymbol{y}]$ if and only if $\boldsymbol{x} \prec \boldsymbol{y}$ for all $[\boldsymbol{x}],[\boldsymbol{y}] \in \frac{\mathcal{X}}{\sim}$. We have that $\preceq$ is an order relation on $\frac{\mathcal{X}}{\sim}$, that is, $(\frac{\mathcal{X}}{\sim},\preceq)$ is a partially ordered set (poset). It is said that $\mathcal{C}$ is a $\preceq$-{\it coverage} if and only if (i) $\mathcal{C} \subseteq \frac{\mathcal{X}}{\sim}$ and (ii) for all $[\boldsymbol{x}] \in \frac{\mathcal{X}}{\sim}$ there exists $[\boldsymbol{y}] \in \mathcal{C}$ such that $[\boldsymbol{y}] \preceq [\boldsymbol{x}]$. Let $[\boldsymbol{y}] \in \frac{\mathcal{X}}{\sim}$, it is said that $[\boldsymbol{y}]$ is {\it minimal} if and only if $\lbrace [\boldsymbol{x}]:[\boldsymbol{x}] \in \frac{\mathcal{X}}{\sim},\;[\boldsymbol{x}] \preceq [\boldsymbol{y}]\rbrace = \lbrace [\boldsymbol{y}] \rbrace$. Let $\mathcal{M}$ be the set of minimal elements for $(\frac{\mathcal{X}}{\sim},\preceq)$, then $\mathcal{M}$ is the (unique) minimum-cardinality $\preceq$-coverage.\\

It is said that $\mathcal{C}$ is a $\prec$-{\it coverage} if and only if (i) $\mathcal{C} \subseteq \mathcal{X}$ and (ii) for all $\boldsymbol{x} \in \mathcal{X}$ there exists $\boldsymbol{y} \in \mathcal{C}$ such that $\boldsymbol{y} \prec \boldsymbol{x}$. Note from \textbf{Remark} \ref{remark:remarkbasica} that $\mathcal{C}$ is a $\prec$-coverage if and only if (i) $\mathcal{C} \subseteq \mathcal{X}$ and (ii) for all $\boldsymbol{x} \in \mathcal{X}$ there exists $\boldsymbol{y} \in \mathcal{C}$ such that $(\boldsymbol{a}^T\boldsymbol{y}, \boldsymbol{c}^T\boldsymbol{y}) \leq (\boldsymbol{a}^T\boldsymbol{x}, \boldsymbol{c}^T\boldsymbol{x})$ for all $(\boldsymbol{a},\boldsymbol{c}) \in \mathcal{U}$. If $\mathcal{C}$ is a $\prec$-coverage let $\delta(\boldsymbol{y}) = |\lbrace \boldsymbol{x}:\boldsymbol{x} \in \mathcal{C},\boldsymbol{x} \prec\boldsymbol{y} \rbrace|$.\\

We have that:
\begin{itemize}
\item  $\mathcal{C}$ is a $\prec$-coverage if and only if there exists $\mathcal{C}^0 \subseteq \mathcal{C}$ such that (i) $\mathcal{M} = \lbrace [\boldsymbol{y}]:\boldsymbol{y} \in \mathcal{C}^0\rbrace$ and (ii) $|\mathcal{C}^0| = |\mathcal{M}|$.
\item $\mathcal{C}$ is a minimum-cardinality $\prec$-coverage if and only if (i) $\mathcal{M} = \lbrace [\boldsymbol{y}]:\boldsymbol{y} \in \mathcal{C}\rbrace$ and (ii) $|\mathcal{C}| = |\mathcal{M}|$. 
\item Let $\mathcal{C}$ be a $\prec$-coverage and let us suppose that $\delta(\boldsymbol{y}) = 1$ for all $\boldsymbol{y} \in \mathcal{C}$ then:  $\mathcal{C}$ is a minimum-cardinality $\prec$-coverage.\\
\end{itemize}

It is said that $\mathcal{C}$ is a $\prec_\epsilon$-{\it coverage} if and only if (i) $\mathcal{C} \subseteq \mathcal{X}$ and (ii) for all $\boldsymbol{x} \in \mathcal{X}$ there exists $\boldsymbol{y} \in \mathcal{C}$ such that $\boldsymbol{y} \prec_\epsilon \boldsymbol{x}$. Note from \textbf{Remark} \ref{remark:remarkbasica} that $\mathcal{C}$ is a $\prec_\epsilon$-coverage if and only if (i) $\mathcal{C} \subseteq \mathcal{X}$ and (ii) for all $\boldsymbol{x} \in \mathcal{X}$ there exists $\boldsymbol{y} \in \mathcal{C}$ such that $(\boldsymbol{a}^T\boldsymbol{y}, \boldsymbol{c}^T\boldsymbol{y}) \leq (\boldsymbol{a}^T\boldsymbol{x}+\epsilon, \boldsymbol{c}^T\boldsymbol{x}+\epsilon)$ for all $(\boldsymbol{a},\boldsymbol{c}) \in \mathcal{U}$.

\begin{remark} \label{remark:UcoverR}
	$\mathcal{C}$ is a $\prec$-coverage if and only if $\mathcal{C}$ is a $\mathcal{U}$-coverage.
\end{remark}

\textbf{Proof}:
\begin{enumerate}
	\item Let us suppose that $\mathcal{C}$ is a $\prec$-coverage. Let $((\boldsymbol{a}^*,\boldsymbol{c}^*),\boldsymbol{x}) \in \mathcal{U} \times \mathcal{X}$. Since $\mathcal{C}$ is a $\prec$-coverage, there exists $\boldsymbol{y} \in \mathcal{C}$ such that $(\boldsymbol{a}^T\boldsymbol{y}, \boldsymbol{c}^T\boldsymbol{y}) \leq (\boldsymbol{a}^T\boldsymbol{x}, \boldsymbol{c}^T\boldsymbol{x})$ for all $(\boldsymbol{a},\boldsymbol{c}) \in \mathcal{U}$, and thus $({\boldsymbol{a}^*}^T\boldsymbol{y},{\boldsymbol{c}^*}^T\boldsymbol{y}) \leq ({\boldsymbol{a}^*}^T\boldsymbol{x},{\boldsymbol{c}^*}^T\boldsymbol{x})$. Therefore, $\mathcal{C}$ is a $\mathcal{U}$-coverage.
	
	\item Let us suppose that $\mathcal{C}$ is a $\mathcal{U}$-coverage and let $\boldsymbol{x} \in \mathcal{X}$. Note that $((\boldsymbol{a}^+(\boldsymbol{x}),\boldsymbol{c}^+(\boldsymbol{x})),\boldsymbol{x}) \in \mathcal{U} \times \mathcal{X}$. Since $\mathcal{C}$ is a $\mathcal{U}$-coverage, there exists $\boldsymbol{y} \in \mathcal{C}$ such that $\boldsymbol{a}^+(\boldsymbol{x})^T\boldsymbol{y} \leq \boldsymbol{a}^+(\boldsymbol{x})^T\boldsymbol{x} = \underline{\mathbf{a}}^T \boldsymbol{x}$ and $\boldsymbol{c}^+(\boldsymbol{x})^T\boldsymbol{y} \leq \boldsymbol{c}^+(\boldsymbol{x})^T\boldsymbol{x} = \underline{\mathbf{c}}^T \boldsymbol{x}$, which implies $\boldsymbol{y} \prec \boldsymbol{x}$. Therefore, $\mathcal{C}$ is a $\prec$-coverage $\blacksquare$
\end{enumerate}

From \textbf{Remark \ref{remark:UcoverR}}, we have that in order to find a minimum-cardinality $\mathcal{U}$-coverage it suffices to find $\mathcal{M}$.

\begin{remark} \label{remark:UcoverReps}
	$\mathcal{C}$ is a $\prec_\epsilon$-coverage if and only if $\mathcal{C}$ is an $\epsilon$-$\mathcal{U}$-coverage.
\end{remark}

\textbf{Proof}:
\begin{enumerate}
	\item Let us suppose that $\mathcal{C}$ is a $\prec_\epsilon$-coverage and let $((\boldsymbol{a}^*,\boldsymbol{c}^*),\boldsymbol{x}) \in \mathcal{U} \times \mathcal{X}$. Since $\mathcal{C}$ is a $\prec_\epsilon$-coverage, there exists $\boldsymbol{y} \in \mathcal{C}$ such that $(\boldsymbol{a}^T\boldsymbol{y}, \boldsymbol{c}^T\boldsymbol{y}) \leq (\boldsymbol{a}^T\boldsymbol{x}+\epsilon, \boldsymbol{c}^T\boldsymbol{x}+\epsilon)$ for all $(\boldsymbol{a},\boldsymbol{c}) \in \mathcal{U}$, and in particular $({\boldsymbol{a}^*}^T\boldsymbol{y},{\boldsymbol{c}^*}^T\boldsymbol{y}) \leq ({\boldsymbol{a}^*}^T\boldsymbol{x}+\epsilon,{\boldsymbol{c}^*}^T\boldsymbol{x}+\epsilon)$. Therefore, $\mathcal{C}$ is an $\epsilon$-$\mathcal{U}$-coverage.
	\item Let us suppose that $\mathcal{C}$ is an $\epsilon$-$\mathcal{U}$-coverage and let $\boldsymbol{x} \in \mathcal{X}$. Note that $((\boldsymbol{a}^+(\boldsymbol{x}),\boldsymbol{c}^+(\boldsymbol{x})),\boldsymbol{x}) \in \mathcal{U} \times \mathcal{X}$. Since $\mathcal{C}$ is an $\epsilon$-$\mathcal{U}$-coverage, there exists $\boldsymbol{y} \in \mathcal{C}$ such that $\boldsymbol{a}^+(\boldsymbol{x})^T\boldsymbol{y} \leq \boldsymbol{a}^+(\boldsymbol{x})^T\boldsymbol{x}+\epsilon = \underline{\mathbf{a}}^T \boldsymbol{x}+\epsilon$ and $\boldsymbol{c}^+(\boldsymbol{x})^T\boldsymbol{y} \leq \boldsymbol{c}^+(\boldsymbol{x})^T\boldsymbol{x}+\epsilon = \underline{\mathbf{c}}^T \boldsymbol{x}+\epsilon$, that is, $\boldsymbol{y} \prec_\epsilon \boldsymbol{x}$. Therefore, $\mathcal{C}$ is a $\prec_\epsilon$-coverage. $\blacksquare$
\end{enumerate}

\begin{remark}\label{remark:coverac}
	If $\mathcal{C}$ is a $\prec$-coverage, then for all $(\boldsymbol{a},\boldsymbol{c}) \in \mathcal{U}$ there exists $\mathcal{C}(\boldsymbol{a},\boldsymbol{c})$, a minimum-cardinality $(\boldsymbol{a},\boldsymbol{c})$-coverage, with $\mathcal{C}(\boldsymbol{a},\boldsymbol{c}) \subseteq \mathcal{C}$.
\end{remark}

\textbf{Proof}: Let $(\boldsymbol{a},\boldsymbol{c}) \in \mathcal{U}$ and let $\mathcal{C}^*(\boldsymbol{a},\boldsymbol{c})$ be a minimum-cardinality $(\boldsymbol{a},\boldsymbol{c})$-coverage. Since $\mathcal{C}$ is a $\prec$-coverage, for all $\boldsymbol{x} \in \mathcal{C}^*(\boldsymbol{a},\boldsymbol{c})$ there exists $\boldsymbol{y} \in \mathcal{C}$ such that $\boldsymbol{y} \prec \boldsymbol{x}$. Then $\mathcal{C}(\boldsymbol{a},\boldsymbol{c}) = \lbrace \boldsymbol{y}: \boldsymbol{y} \in \mathcal{C},\boldsymbol{y} \prec \boldsymbol{x},\boldsymbol{x} \in \mathcal{C}^*(\boldsymbol{a},\boldsymbol{c}) \rbrace$ is a $(\boldsymbol{a},\boldsymbol{c})$-coverage with $|\mathcal{C}(\boldsymbol{a},\boldsymbol{c})| \leq |\mathcal{C}^*(\boldsymbol{a},\boldsymbol{c})|\;\blacksquare$\\

Note that if $\mathcal{C}$ is a $\prec$-coverage then  a minimum-cardinality $(\boldsymbol{a},\boldsymbol{c})$-coverage included in $\mathcal{C}$ may be found easily. 

\begin{remark}\label{remark:poset}
	Let $J^+ = \lbrace j \in [n]:\mbox{ either } \underline{\mathbf{a}}_j < \overline{\mathbf{a}}_j \mbox{ or } \underline{\mathbf{c}}_j < \overline{\mathbf{c}}_j \rbrace$. If $J^+=[n]$, then $(\mathcal{X},\prec)$ is a poset.
\end{remark}

\textbf{Proof}: Let us suppose that $J^+=[n]$. Let $\boldsymbol{x},\boldsymbol{y} \in \mathcal{X}$ and suppose that $\boldsymbol{x} \sim \boldsymbol{y}$. By definition of $\sim$, this implies $\boldsymbol{x} \prec \boldsymbol{y}$ and $\boldsymbol{y} \prec \boldsymbol{x}$, meaning $\boldsymbol{a}^+(\boldsymbol{y})^T\boldsymbol{x} \leq \underline{\mathbf{a}}^T\boldsymbol{y}$, $\boldsymbol{c}^+(\boldsymbol{y})^T\boldsymbol{x} \leq \underline{\mathbf{c}}^T\boldsymbol{y}$, $\boldsymbol{a}^+(\boldsymbol{x})^T\boldsymbol{y} \leq \underline{\mathbf{a}}^T\boldsymbol{x}$, and $\boldsymbol{c}^+(\boldsymbol{x})^T\boldsymbol{y} \leq \underline{\mathbf{c}}^T\boldsymbol{x}$. Since $J^+ = [n]$, by applying \textbf{Remark \ref{remark:remarkJxy}} to both objective bounds across all coordinates, we obtain $\boldsymbol{x} = \boldsymbol{y}$. Hence, each equivalence class is a singleton $[\boldsymbol{x}] = \lbrace \boldsymbol{x} \rbrace$ and then $(\mathcal{X},\prec)$ is a poset $\blacksquare$

\section{Algorithms}

\subsection{Algorithms to find minimal elements}

\begin{remark}\label{remark:problemR}  
	Let $J^+ = \lbrace j \in [n]: \mbox{ either }\underline{\mathbf{a}}_j  < \overline{\mathbf{a}}_j \mbox{ or } \underline{\mathbf{c}}_j < \overline{\mathbf{c}}_j \rbrace$, let $\boldsymbol{x} \in \mathcal{X}$, let $K^1(\boldsymbol{x})=\lbrace j:j \in J^+,\boldsymbol{x}_j = 1 \rbrace$, let $K^0(\boldsymbol{x})=\lbrace j:j \in J^+,\boldsymbol{x}_j = 0 \rbrace$ and let $\mathcal{R}(\boldsymbol{x})$ be a problem in $\boldsymbol{y}$ defined as follows:
	\begin{flalign*}
		\min&\;\sum_{j \in K^1(\boldsymbol{x})}\boldsymbol{y}_j - \sum_{j \in K^0(\boldsymbol{x})}\boldsymbol{y}_j &\;& &\mathcal{R}(\boldsymbol{x})&\\
		s.t.&\;\;\boldsymbol{y} \prec \boldsymbol{x}&\\
		&\;\;\boldsymbol{y} \in \mathcal{X}&\
	\end{flalign*}
	and let $\boldsymbol{y}^*$ be an optimal solution then:
	\begin{enumerate}
		\item $[\boldsymbol{y}^*] \preceq [\boldsymbol{x}]$
		\item If $v(\mathcal{R}(\boldsymbol{x})) < |K^1(\boldsymbol{x})|$ then  $[\boldsymbol{y}^*] \neq [\boldsymbol{x}]$
		\item If $v(\mathcal{R}(\boldsymbol{x})) = |K^1(\boldsymbol{x})|$ then $\boldsymbol{y}_j = \boldsymbol{x}_j$ for all $j \in J^+$ for all $\boldsymbol{y} \in \mathcal{X}$ such that $\boldsymbol{y} \prec \boldsymbol{x}$ 
		\item If $v(\mathcal{R}(\boldsymbol{x})) = |K^1(\boldsymbol{x})|$ and $J^+ = [n]$ then $\boldsymbol{x} \in \mathcal{M}$
		
	\end{enumerate}	
\end{remark}
\textbf{Proof}:
\begin{enumerate}
	\item From the definition of $\mathcal{R}(\boldsymbol{x})$.
	\item If $[\boldsymbol{y}^*] = [\boldsymbol{x}]$ then we have that $\boldsymbol{y}^* \prec \boldsymbol{x}$ and $\boldsymbol{x} \prec \boldsymbol{y}^*$ and from \textbf{Remark} \ref{remark:remarkJxy} we have that $\boldsymbol{x}_j = \boldsymbol{y}^*_j$ for all $j \in J^+$ therefore $v(\mathcal{R}(\boldsymbol{x})) = \sum_{j \in K^1(\boldsymbol{x})}\boldsymbol{y}^*_j - \sum_{j \in K^0(\boldsymbol{x})}\boldsymbol{y}^*_j = |K^1(\boldsymbol{x})|$.
	\item   Let $\boldsymbol{y} \in  F(\mathcal{R}(\boldsymbol{x}))$, hence $\sum_{j \in K^1(\boldsymbol{x})}\boldsymbol{y}_j - \sum_{j \in K^0(\boldsymbol{x})}\boldsymbol{y}_j \geq v(\mathcal{R}(\boldsymbol{x})) = |K^1(\boldsymbol{x})|$, therefore $\boldsymbol{y}_j = \boldsymbol{x}_j$ for all $j \in J^+$.
	\item Let us suppose that $v(\mathcal{R}(\boldsymbol{x})) = |K^1(\boldsymbol{x})|$ and $J^+ = [n]$ then from \textbf{Remark} \ref{remark:poset} we have that $(\mathcal{X},\prec)$ is a poset and from 3 we have that $\boldsymbol{x} = \boldsymbol{y}$ for all $\boldsymbol{y} \in \mathcal{X}$ such that $\boldsymbol{y} \prec \boldsymbol{x}$, therefore $\boldsymbol{x}$ is a minimal element $\blacksquare$   
\end{enumerate}

In order to solve $\mathcal{R}(\boldsymbol{x})$ with a solver (for example with CPLEX) note that  $\mathcal{R}(\boldsymbol{x})$ may be rewritten as follows:

\begin{flalign*}
	\min&\;\sum_{j \in K^1(\boldsymbol{x})}\boldsymbol{y}_j - \sum_{j \in K^0(\boldsymbol{x})}\boldsymbol{y}_j &\;& &\mathcal{R}(\boldsymbol{x})&\\
	s.t.&\;\;\boldsymbol{a}^+(\boldsymbol{x})^T\boldsymbol{y} \leq \underline{\mathbf{a}}^T \boldsymbol{x}\\
	&\;\;\boldsymbol{c}^+(\boldsymbol{x})^T\boldsymbol{y} \leq \underline{\mathbf{c}}^T \boldsymbol{x}\\
	&\;\;\boldsymbol{y} \in \mathcal{X}&\
\end{flalign*}

and since $\boldsymbol{x}$ is fixed in problem $\mathcal{R}(\boldsymbol{x})$ we have that $\mathcal{R}(\boldsymbol{x})$ is a 0-1-MILP problem.\\

Below we present an algorithm to find a minimal element starting from $\boldsymbol{x}^0 \in \mathcal{X}$ when $J^+=[n]$.    

\begin{algorithm}[H]
	\caption{$minimal(\boldsymbol{x}^0)/J^+=[n]$}
	\hypertarget{target:MinimalR}{}
	\renewcommand{\thealgorithm}{}
	\begin{algorithmic}[1]
		\State $stop \leftarrow 0$
		\State $\boldsymbol{x}  \leftarrow \boldsymbol{x}^0$
		\While{$stop=0$}
		\State  solve $\mathcal{R}(\boldsymbol{x})$
		\If {$v(\mathcal{R}(\boldsymbol{x})) = |K^1(\boldsymbol{x})|$} 
		\State $stop \leftarrow 1$
		\Else
		\State let $\boldsymbol{y}$ be an optimal solution for $\mathcal{R}(\boldsymbol{x})$,
		\State {$\boldsymbol{x} \leftarrow \boldsymbol{y}$}
		\EndIf
		\EndWhile
		\State \textbf{return} $\boldsymbol{x}$
	\end{algorithmic}
\end{algorithm}

\begin{remark}\label{remark:minimal1}
	Let $\boldsymbol{x}^0 \in \mathcal{X}$, let $J^+=[n]$ and let $\boldsymbol{x}$ be the output of \hyperlink{target:MinimalR}{$minimal(\boldsymbol{x}^0)/J^+=[n]$}, then $\boldsymbol{x} \in \mathcal{M}$.
\end{remark}
\textbf{Proof}:
Let $\mathcal{Y}$ be the set of solutions generated by the algorithm before step 4. If 
$v(\mathcal{R}(\boldsymbol{x})) < |K^1(\boldsymbol{x})|$ then we have $\boldsymbol{y} \notin \mathcal{Y}$ for the new generated solution. Since $\mathcal{X}$ is a finite set any sequence generated by the algorithm must be finite. The loop terminates when $v(\mathcal{R}(\boldsymbol{x})) = |K^1(\boldsymbol{x})|$ and from \textbf{Remark \ref{remark:problemR}} we have that $\boldsymbol{x} \in \mathcal{M}\;\blacksquare$.

\begin{remark}\label{remark:problemS}
	Let $\boldsymbol{x} \in \mathcal{X}$ and let $\mathcal{S}(\boldsymbol{x})$ be a problem in $\boldsymbol{y}$ defined as follows:
	\begin{flalign*}
		find&\;\;\boldsymbol{y}&\;& &\mathcal{S}(\boldsymbol{x})&\\
		s.t.&\;\;\boldsymbol{y} \prec \boldsymbol{x}&\\
		&\;\;\boldsymbol{x} \nprec_+ \boldsymbol{y} \mbox{ or } \boldsymbol{x} \nprec^+ \boldsymbol{y}&\\
		&\;\;\boldsymbol{y} \in \mathcal{X}&\
	\end{flalign*}
	then:
	\begin{enumerate}
		\item If $F(\mathcal{S}(\boldsymbol{x})) = \emptyset$ then $[\boldsymbol{x}] \in \mathcal{M}$
		\item Let $\boldsymbol{y}^* \in F(\mathcal{S}(\boldsymbol{x}))$  then $[\boldsymbol{y}^*] \preceq [\boldsymbol{x}]$ and $[\boldsymbol{y}^*] \neq [\boldsymbol{x}]$ 
	\end{enumerate}
\end{remark}
\textbf{Proof}:
\begin{enumerate}
	\item Note that if $\boldsymbol{y} \in F(\mathcal{S}(\boldsymbol{x}))$ then $\boldsymbol{y} \prec \boldsymbol{x}$ and $\boldsymbol{x} \nprec\boldsymbol{y}$ then if $F(\mathcal{S}(\boldsymbol{x})) = \emptyset$ we have that $\boldsymbol{x} \prec \boldsymbol{y}$ for all $\boldsymbol{y} \in \mathcal{X}$ such that $\boldsymbol{y} \prec \boldsymbol{x}$, hence $[\boldsymbol{x}] \in \mathcal{M}$.
	\item If $\boldsymbol{y}^* \in F(\mathcal{S}(\boldsymbol{x}))$ then $\boldsymbol{y}^* \prec \boldsymbol{x}$, hence $[\boldsymbol{y}^*] \preceq [\boldsymbol{x}]$. Since $\boldsymbol{x}\nprec\boldsymbol{y}^*$ we have that $[\boldsymbol{y}^*] \neq [\boldsymbol{x}]\;\blacksquare$
\end{enumerate}

In order to solve $\mathcal{S}(\boldsymbol{x})$ with a solver (for example with CPLEX) note that  $\mathcal{S}(\boldsymbol{x})$ may be rewritten as follows:

\begin{flalign*}
	find&\;\;\boldsymbol{y}&\;& &\mathcal{S}(\boldsymbol{x})&\\
	s.t.&\;\;\boldsymbol{a}^+(\boldsymbol{x})^T\boldsymbol{y} \leq \underline{\mathbf{a}}^T \boldsymbol{x}\\
	&\;\;\boldsymbol{c}^+(\boldsymbol{x})^T\boldsymbol{y} \leq \underline{\mathbf{c}}^T \boldsymbol{x}\\
	&\;\;\boldsymbol{a}^+(\boldsymbol{y})^T\boldsymbol{x} \geq \underline{\mathbf{a}}^T \boldsymbol{y}+1&\mbox{if }\boldsymbol{z}=1\\ &\;\;\boldsymbol{c}^+(\boldsymbol{y})^T\boldsymbol{x} \geq \underline{\mathbf{c}}^T \boldsymbol{y}+1&\mbox{if }\boldsymbol{z}=0\\
	&\;\;\boldsymbol{y} \in \mathcal{X},\;\boldsymbol{z} \in \lbrace 0,1 \rbrace&\
\end{flalign*}

and since $\boldsymbol{x}$ is fixed in problem $\mathcal{S}(\boldsymbol{x})$ we have that $\mathcal{S}(\boldsymbol{x})$ is a disjunctive 0-1-MILP problem.\\

Below we present an algorithm to find a minimal element starting from $\boldsymbol{x}^0 \in \mathcal{X}$.

\begin{algorithm}[H]
	\caption{$minimal(\boldsymbol{x}^0)$}
	\hypertarget{target:MinimalS}{}
	\renewcommand{\thealgorithm}{}
	\begin{algorithmic}[1]
		\State $stop \leftarrow 0$
		\State $\boldsymbol{x}  \leftarrow \boldsymbol{x}^0$		
		\While{$stop=0$}
		\State  solve $\mathcal{S}(\boldsymbol{x})$
		\If {$F(\mathcal{S}(\boldsymbol{x})) = \emptyset$} 
		\State $stop \leftarrow 1$
		\Else
		\State let $\boldsymbol{y} \in F(\mathcal{S}(\boldsymbol{x}))$
		\State {$\boldsymbol{x} \leftarrow \boldsymbol{y}$}
		\EndIf
		\EndWhile
		\State \textbf{return} $\boldsymbol{x}$
	\end{algorithmic}
\end{algorithm}

\begin{remark}\label{remark:minimal2}
	Let $\boldsymbol{x}^0 \in \mathcal{X}$ and let $\boldsymbol{x}$ be the output of \hyperlink{target:MinimalS}{$minimal(\boldsymbol{x}^0)$}, then $[\boldsymbol{x}] \in \mathcal{M}$.
\end{remark}

\textbf{Proof}:
Let $\mathcal{Y}$ be the set of solutions generated by the algorithm before step 4. If 
$F(\mathcal{S}(\boldsymbol{x})) \neq \emptyset$ then we have $\boldsymbol{y} \notin \mathcal{Y}$ for the new generated solution. Since $\mathcal{X}$ is a finite set any sequence generated by the algorithm must be finite. The loop terminates when  $F(\mathcal{S}(\boldsymbol{x})) = \emptyset$ and from \textbf{Remark \ref{remark:problemS}} we have that $[\boldsymbol{x}] \in \mathcal{M}\;\blacksquare$.

\subsection{Disjunctive algorithm}

Let $\mathcal{C} \subseteq \mathcal{X}$ such that $\lbrace [\boldsymbol{y}]:\boldsymbol{y} \in \mathcal{C}\rbrace \subseteq \mathcal{M}$, let $\epsilon \in \mathbb{N}\cup\lbrace 0 \rbrace$ and let $\mathcal{Q}(\mathcal{C},\epsilon)$ be a problem in $(\boldsymbol{x},\boldsymbol{z})$  defined as follows:
\begin{flalign*}
	find&\;\;\boldsymbol{x} &\;&&\mathcal{Q}(\mathcal{C},\epsilon)&\\
	s.t.&\;\;\boldsymbol{y} \nprec^\epsilon_+ \boldsymbol{x}&\mbox{if }\boldsymbol{z}(\boldsymbol{y})=1\;\forall \boldsymbol{y} \in \mathcal{C}\\ &\;\;\boldsymbol{y} \nprec^+_\epsilon \boldsymbol{x}&\mbox{if }\boldsymbol{z}(\boldsymbol{y})=0\;\forall \boldsymbol{y} \in \mathcal{C}&\\
	&\;\;\boldsymbol{x} \in \mathcal{X},\;\boldsymbol{z}(\boldsymbol{y}) \in \lbrace 0,1 \rbrace\;\forall \boldsymbol{y} \in \mathcal{C}&\\
\end{flalign*}
\begin{remark}\label{remark:problemQeps}
	$\;$
	\begin{enumerate}
	\item If $\epsilon \in \mathbb{N} \cup \lbrace 0 \rbrace$ and $F(\mathcal{Q}(\mathcal{C},\epsilon))  = \emptyset$ then $\mathcal{C}$ is an  $\epsilon$-$\mathcal{U}$-coverage.
	\item If $\epsilon = 0$ and $F(\mathcal{Q}(\mathcal{C},0))  = \emptyset$ then $\lbrace [\boldsymbol{y}]:\boldsymbol{y} \in \mathcal{C}\rbrace =  \mathcal{M}$
	\item Let $\boldsymbol{x} \in F(\mathcal{Q}(\mathcal{C},\epsilon))$ and let $[\boldsymbol{y}^*] \in \mathcal{M}$ such that $[\boldsymbol{y}^*] \preceq [\boldsymbol{x}]$ then $[\boldsymbol{y}^*] \notin \lbrace [\boldsymbol{y}]:\boldsymbol{y} \in \mathcal{C}\rbrace$
	\end{enumerate} 
\end{remark}

\textbf{Proof}: 
\begin{enumerate}
\item Let us suppose that $\epsilon \in \mathbb{N} \cup \lbrace 0 \rbrace$. If $F(\mathcal{Q}(\mathcal{C},\epsilon))  = \emptyset$ then for all $\boldsymbol{x} \in \mathcal{X}$ there exists $\boldsymbol{y} \in \mathcal{C}$ such that
$\boldsymbol{a}^+(\boldsymbol{x})^T\boldsymbol{y} \leq \underline{\mathbf{a}}^T \boldsymbol{x}+\epsilon$ and $\boldsymbol{c}^+(\boldsymbol{x})^T\boldsymbol{y} \leq \underline{\mathbf{c}}^T \boldsymbol{x}+\epsilon$, therefore  $\mathcal{C}$ is an   $\epsilon$-$\mathcal{U}$-coverage .
\item Let us suppose that $\epsilon = 0$. If $F(\mathcal{Q}(\mathcal{C},0))  = \emptyset$ we have that for all $\boldsymbol{x} \in \mathcal{X}$ there exists $\boldsymbol{y} \in \mathcal{C}$ such that $\boldsymbol{y} \prec \boldsymbol{x}$, therefore $\mathcal{C}$ is a $\prec$-coverage and since $\lbrace [\boldsymbol{y}]:\boldsymbol{y} \in \mathcal{C}\rbrace \subseteq \mathcal{M}$ we have that $\lbrace [\boldsymbol{y}]:\boldsymbol{y} \in \mathcal{C}\rbrace = \mathcal{M}$.
\item If $[\boldsymbol{y}^*] \in \lbrace [\boldsymbol{y}]:\boldsymbol{y} \in \mathcal{C}\rbrace$ then there exists $\boldsymbol{y} \in \mathcal{C}$ such that $\boldsymbol{y} \prec \boldsymbol{y}^* \prec \boldsymbol{x}$ and then $\boldsymbol{y} \prec_\epsilon \boldsymbol{x}$ and $\boldsymbol{x} \notin F(\mathcal{Q}(\mathcal{C},\epsilon))\;\blacksquare$
 \end{enumerate}

In order to solve $\mathcal{Q}(\mathcal{C},\epsilon)$ with a solver (for example with CPLEX) note that $\mathcal{Q}(\mathcal{C},\epsilon)$ may be rewritten as follows:

\begin{flalign*}
	find&\;\;\boldsymbol{x} &\;&&\mathcal{Q}(\mathcal{C},\epsilon)&\\
	s.t.&\;\;\boldsymbol{a}^+(\boldsymbol{x})^T\boldsymbol{y} \geq \underline{\mathbf{a}}^T \boldsymbol{x}+\epsilon+1&\mbox{if }\boldsymbol{z}(\boldsymbol{y})=1\;\forall \boldsymbol{y} \in \mathcal{C}\\ &\;\;\boldsymbol{c}^+(\boldsymbol{x})^T\boldsymbol{y} \geq \underline{\mathbf{c}}^T \boldsymbol{x}+\epsilon+1&\mbox{if }\boldsymbol{z}(\boldsymbol{y})=0\;\forall \boldsymbol{y} \in \mathcal{C}&\\
	&\;\;\boldsymbol{x} \in \mathcal{X},\;\boldsymbol{z}(\boldsymbol{y}) \in \lbrace 0,1 \rbrace\;\forall \boldsymbol{y} \in \mathcal{C}&\\
\end{flalign*}

and since $\mathcal{C}$ is fixed in problem $\mathcal{Q}(\mathcal{C},\epsilon)$ we have that $\mathcal{Q}(\mathcal{C},\epsilon)$ is a disjunctive 0-1-MILP problem.\\

Below we present an algorithm to find, given $\epsilon \in \mathbb{N}\cup \lbrace 0 \rbrace$, an $\epsilon$-$\mathcal{U}$-coverage. If $\epsilon=0$ the algorithm finds a minimum-cardinality $\mathcal{U}$-coverage. The algorithm finds in each iteration a solution to be included in the coverage and adds a new binary variable and two constraints to be used in the next iterations. If $J^+ = [n]$ then \hyperlink{target:MinimalR}{$minimal(\boldsymbol{x})/J^+=[n]$} may be used instead of \hyperlink{target:MinimalS}{$minimal(\boldsymbol{x})$}. 

\begin{algorithm}[H]
	\caption{$Disj(\epsilon)$}
	\hypertarget{target:Disjunctive}{}
	\renewcommand{\thealgorithm}{}
	\begin{algorithmic}[1]
		\State Let $\boldsymbol{x} \in \mathcal{X}$
		\State Use \hyperlink{target:MinimalS}{$minimal(\boldsymbol{x})$} to obtain $\boldsymbol{y} \in \mathcal{X}$ such that $\boldsymbol{y} \prec \boldsymbol{x}$ and $[\boldsymbol{y}] \in \mathcal{M}$
		\State $\mathcal{C} \leftarrow \lbrace \boldsymbol{y} \rbrace$ 
		\State $stop \leftarrow 0$
		\While{$stop=0$}
		\State  solve $\mathcal{Q}(\mathcal{C},\epsilon)$
		\If {$F(\mathcal{Q}(\mathcal{C},\epsilon)) = \emptyset$} 
		\State $stop \leftarrow 1$
		\Else
		\State Let $\boldsymbol{x}$ be an optimal solution for $\mathcal{Q}(\mathcal{C},\epsilon)$,
		\State use \hyperlink{target:MinimalS}{$minimal(\boldsymbol{x})$} to obtain $\boldsymbol{y} \in \mathcal{X}$ such that $\boldsymbol{y} \prec \boldsymbol{x}$ and $[\boldsymbol{y}] \in \mathcal{M}$
		\State $\mathcal{C} \leftarrow \mathcal{C} \cup \lbrace\boldsymbol{y} \rbrace$ 
		\EndIf
		\EndWhile
		\State \textbf{return} $\mathcal{C}$
	\end{algorithmic}
\end{algorithm}

\begin{remark}\label{remark:algorithmDisjunction}
	$\;$
	\begin{enumerate}
		\item \hyperlink{target:Disjunctive}{$Disj(\epsilon)$} is finite.
		\item Let $\mathcal{C}$ be the output of \hyperlink{target:Disjunctive}{$Disj(\epsilon)$}	then:

	\begin{enumerate}
		\item If $\epsilon=0$ then $\mathcal{C}$ is a minimum-cardinality  $\mathcal{U}$-coverage
		\item If $\epsilon \in \mathbb{N} \cup \lbrace 0 \rbrace$ then $\mathcal{C}$ is an $\epsilon$-$\mathcal{U}$-coverage. 
	\end{enumerate}
\end{enumerate}
\end{remark}

\textbf{Proof}:
\begin{enumerate}
	\item  If $F(\mathcal{Q}(\mathcal{C},\epsilon)) \neq \emptyset$ then from \textbf{Remark \ref{remark:problemQeps}} we have that $[\boldsymbol{y}] \notin \{[\boldsymbol{w}] : \boldsymbol{w} \in \mathcal{C}\}$ for the new generated solution. Since $\mathcal{X}$ is a finite set the algorithm stops in a finite number of iterations.
	\item Let $\mathcal{C}$ be the output of \hyperlink{target:Disjunctive}{$Disj(\epsilon)$}	then:
\begin{enumerate}
\item Let us suppose that $\epsilon = 0$. By definition of the algorithm we have that $\{[\boldsymbol{y}] : \boldsymbol{y} \in \mathcal{C}\} \subseteq \mathcal{M}$ before step 6. If $F(\mathcal{Q}(\mathcal{C},0)) = \emptyset$ then from  \textbf{Remark \ref{remark:problemQeps}} we have that $\{[\boldsymbol{y}] : \boldsymbol{y} \in \mathcal{C}\} = \mathcal{M}$ and then $\mathcal{C}$ is a minimum-cardinality $\mathcal{U}$-coverage.  
\item Let us suppose that $\epsilon \in \mathbb{N} \cup \lbrace 0 \rbrace$. By definition of the algorithm we have that $\{[\boldsymbol{y}] : \boldsymbol{y} \in \mathcal{C}\} \subseteq \mathcal{M}$ before step 6. If $F(\mathcal{Q}(\mathcal{C},\epsilon)) = \emptyset$ then from  \textbf{Remark \ref{remark:problemQeps}} we have that  $\mathcal{C}$ is an $\epsilon$-$\mathcal{U}$-coverage $\blacksquare$ 
\end{enumerate}
\end{enumerate}

\subsection{Driven-disjunctive algorithm}

Drawing inspiration from the concepts presented in the work of \cite{Bektas2018}, we have developed an algorithm  to find an $\epsilon$-$\mathcal{U}$-coverage that does not incorporate any additional variables. If $\epsilon=0$ the algorithm finds a minimum-cardinality $\mathcal{U}$-coverage.\\

Let $\mathcal{C} \subseteq \mathcal{X}$ such that $\lbrace [\boldsymbol{y}]:\boldsymbol{y} \in \mathcal{C}\rbrace \subseteq \mathcal{M}$, let $\epsilon \in \mathbb{N} \cup \lbrace 0 \rbrace $, let $\boldsymbol{z}:\mathcal{C} \rightarrow {\lbrace 0,1 \rbrace}^{|\mathcal{C}|}$ and let $\mathcal{T}(\mathcal{C},\epsilon,\boldsymbol{z})$ be a problem in $\boldsymbol{x}$ defined as follows:
\begin{flalign*}
	find&\;\;\boldsymbol{x} &\;&&\mathcal{T}(\mathcal{C},\epsilon,\boldsymbol{z})&\\
	s.t.&\;\;\boldsymbol{y} \nprec^\epsilon_+ \boldsymbol{x}&\mbox{ if }\boldsymbol{z}(\boldsymbol{y})=1\;\forall \boldsymbol{y} \in \mathcal{C}\\ &\;\;\boldsymbol{y} \nprec^+_\epsilon \boldsymbol{x}&\mbox{ if }\boldsymbol{z}(\boldsymbol{y})=0\;\forall \boldsymbol{y} \in \mathcal{C}&\\
	&\;\;\boldsymbol{x} \in \mathcal{X}&\
\end{flalign*}

In order to solve $\mathcal{T}(\mathcal{C},\epsilon,\boldsymbol{z})$  with a solver (for example with CPLEX) note that $\mathcal{T}(\mathcal{C},\epsilon,\boldsymbol{z})$ may be rewritten as follows:
\begin{flalign*}
	find&\;\;\boldsymbol{x} &\;&&\mathcal{T}(\mathcal{C},\epsilon,\boldsymbol{z})&\\
	s.t.&\;\;\boldsymbol{a}^+(\boldsymbol{x})^T\boldsymbol{y} \geq \underline{\mathbf{a}}^T \boldsymbol{x}+\epsilon+1&\mbox{ if }\boldsymbol{z}(\boldsymbol{y})=1\;\forall \boldsymbol{y} \in \mathcal{C}\\ &\;\;\boldsymbol{c}^+(\boldsymbol{x})^T\boldsymbol{y} \geq \underline{\mathbf{c}}^T \boldsymbol{x}+\epsilon+1&\mbox{ if }\boldsymbol{z}(\boldsymbol{y})=0\;\forall \boldsymbol{y} \in \mathcal{C}&\\
	&\;\;\boldsymbol{x} \in \mathcal{X}&\
\end{flalign*}

and since $\boldsymbol{z}$ and $\mathcal{C}$ are  fixed in problem $\mathcal{T}(\mathcal{C},\epsilon,\boldsymbol{z})$ we have that $\mathcal{T}(\mathcal{C},\epsilon,\boldsymbol{z})$ is a 0-1-MILP problem.\\

Below we present an algorithm to find, given $\epsilon \in \mathbb{N}\cup \lbrace 0 \rbrace$, an $\epsilon$-$\mathcal{U}$-coverage. If $\epsilon=0$ the algorithm finds a minimum-cardinality $\mathcal{U}$-coverage. If $J^+ = [n]$ then \hyperlink{target:MinimalR}{$minimal(\boldsymbol{x})/J^+=[n]$} may be used instead of \hyperlink{target:MinimalS}{$minimal(\boldsymbol{x})$}

\begin{algorithm}[H]
	\caption{$Driven(\epsilon)$}
	\hypertarget{target:Branch and Driven-Disjunction}{}
	\renewcommand{\thealgorithm}{}
	\begin{algorithmic}[1]
		\State Let $\boldsymbol{x} \in \mathcal{X}$
		\State Use \hyperlink{target:MinimalS}{$minimal(\boldsymbol{x})$} to obtain $\boldsymbol{y} \in \mathcal{X}$ such that $\boldsymbol{y} \prec \boldsymbol{x}$ and $[\boldsymbol{y}] \in \mathcal{M}$
		\State $\mathcal{C} \leftarrow \lbrace \boldsymbol{y} \rbrace$, $\mathcal{C}^* \leftarrow \lbrace \boldsymbol{y} \rbrace$ and $\boldsymbol{z}(\boldsymbol{y}) \leftarrow 1$
		\State the {\it brother}-problem (the same $\mathcal{C}$ with $\boldsymbol{z}(\boldsymbol{y})=0$) is included in the $list$
		\State $stop \leftarrow 0$
		\While{$stop=0$}
		\State  solve $\mathcal{T}(\mathcal{C},\epsilon,\boldsymbol{z})$
		\If {$F(\mathcal{T}(\mathcal{C},\epsilon,\boldsymbol{z})) = \emptyset$}
		\State ({\it the node is fathomed}) 
		\If {$list = \emptyset$}
		\State $stop \leftarrow 1$
		\Else 
		\State Choose the next
		$(\mathcal{C},\boldsymbol{z})$ from the $list$
		\State $list \leftarrow list - \lbrace (\mathcal{C},\boldsymbol{z}) \rbrace$ 
		\EndIf
		\Else
		\State Let $\boldsymbol{x} \in F(\mathcal{T}(\mathcal{C},\epsilon,\boldsymbol{z}))$,
		\State Use \hyperlink{target:MinimalS}{$minimal(\boldsymbol{x})$} to obtain $\boldsymbol{y} \in \mathcal{X}$ such that $\boldsymbol{y} \prec \boldsymbol{x}$ and $[\boldsymbol{y}] \in \mathcal{M}$
		\State ({\it branching})
		\State  $\mathcal{C} \leftarrow \mathcal{C} \cup \lbrace\boldsymbol{y} \rbrace$ and $\boldsymbol{z}(\boldsymbol{y}) \leftarrow 1$
		\State the {\it brother}-problem (the same $\mathcal{C}$ with $\boldsymbol{z}(\boldsymbol{y})=0$) is included in the $list$
		\If {$\boldsymbol{y} \nprec \boldsymbol{w}\; \forall \boldsymbol{w} \in \mathcal{C}^*$}
		\State $\mathcal{C}^* \leftarrow \mathcal{C}^* \cup \lbrace \boldsymbol{y} \rbrace$
		\EndIf
		\State Add a cut to $\mathcal{X}$ to exclude $\boldsymbol{y}$ from future subproblems
		\EndIf
		\EndWhile
		\State \textbf{return} $\mathcal{C}^*$
	\end{algorithmic}
\end{algorithm}

\begin{remark}\label{remark:algorithmDrivenDisjunction}
	$\;$
	\begin{enumerate}
	\item \hyperlink{target:Branch and Driven-Disjunction}{$Driven(\epsilon)$}	is finite.
	\item
	Let $\mathcal{C}^*$ be the output of \hyperlink{target:Branch and Driven-Disjunction}{$Driven(\epsilon)$}	then:
	\begin{enumerate}
		\item $\mathcal{C}^*$ is an $\epsilon$-$\mathcal{U}$-coverage. 		
		\item If $\epsilon=0$ then $\mathcal{C}^*$ is a minimum-cardinality $\mathcal{U}$-coverage
	\end{enumerate}
\end{enumerate}
\end{remark}

\textbf{Proof}:
\begin{enumerate} 
	\item  The algorithm generates a binary tree. At each node of the tree (the {\it parent node}), the corresponding problem $\mathcal{T}(\mathcal{C},\epsilon,\boldsymbol{z})$ is solved. If $F(\mathcal{T}(\mathcal{C},\epsilon,\boldsymbol{z})) \neq \emptyset$ and if $\boldsymbol{x}$ is the solution found, a minimal element $\boldsymbol{y}$ is generated by using \hyperlink{target:MinimalS}{$minimal(\boldsymbol{x})$}. The parent node gives rise to two {\it child nodes}, corresponding to the inclusion of the constraints $\boldsymbol{z}(\boldsymbol{y})=1$ for the first child and $\boldsymbol{z}(\boldsymbol{y})=0$ for the second. From the definition of $F(\mathcal{T}(\mathcal{C},\epsilon,\boldsymbol{z}))$ for each node, it follows immediately that the minimal elements defining a given branch cannot be repeated within that same branch. While each minimal element may appear multiple times throughout the tree, it can appear at most once in any specific branch. Given that the set of minimal elements is finite, each branch has a finite length, and eventually, the problem $\mathcal{T}(\mathcal{C},\epsilon,\boldsymbol{z})$ becomes infeasible. Since the tree is binary and each branch has a finite length, the tree itself is finite. It follows that \hyperlink{target:Branch and Driven-Disjunction}{$Driven(\epsilon)$} is finite.\\

	\item Let $\mathcal{C}^*$ be the output of \hyperlink{target:Branch and Driven-Disjunction}{$Driven(\epsilon)$}.  
	\begin{enumerate}
		\item  Note that if $\boldsymbol{y}$ is discarded in steps 22 and 23 and does not become an element of $\mathcal{C}^*$, it is because there exists some $\boldsymbol{w} \in \mathcal{C}^*$ such that $[\boldsymbol{y}] = [\boldsymbol{w}]$ (since both are minimal); therefore, for the purpose of covering an element of $\mathcal{X}$, it makes no difference whether $\mathcal{C}$ or $\mathcal{C}^*$ is used. Suppose that $\mathcal{C}$ does not constitute an $\epsilon$-$\mathcal{U}$-covering; then, there exists an $\boldsymbol{x} \in \mathcal{X}$ such that $\boldsymbol{y} \nprec_\epsilon \boldsymbol{x}$ for all $\boldsymbol{y} \in \mathcal{C}$. In this case, $\boldsymbol{x}$ would have to be feasible in one of the branches of the tree. Given that all branches must have been closed due to infeasibility, we arrive at a contradiction; consequently, $\mathcal{C}$ must be an $\epsilon$-$\mathcal{U}$-covering, and the same holds for $\mathcal{C}^*$.\\

		\item If $\epsilon=0$ we have that $\mathcal{C}^*$ is a $\mathcal{U}$-coverage and by steps 22 and 23 we have that $\delta(\boldsymbol{y}) = 1$ for all $\boldsymbol{y} \in \mathcal{C}^*$ and then $\mathcal{C}^*$ is a minimum-cardinality $\mathcal{U}$-coverage $\blacksquare$   
\end{enumerate}
\end{enumerate}

In order to prevent previously generated minimal elements from reappearing directly as feasible solutions $\boldsymbol{x}$ in subsequent subproblems $\mathcal{T}(\mathcal{C},\epsilon,\boldsymbol{z})$, and thus save memory and computing time, a standard cut excluding $\boldsymbol{y}$ is added to $\mathcal{X}$ in step 25:
\[
\sum_{j \in [n]: \boldsymbol{y}_j = 1} (1 - \boldsymbol{x}_j) + \sum_{j \in [n]: \boldsymbol{y}_j = 0} \boldsymbol{x}_j \geq 1.
\]
Note that while this cut prevents $\boldsymbol{y}$ from being generated directly as a solution $\boldsymbol{x}$, it does not prevent $\boldsymbol{y}$ from reappearing as the minimal element associated with a non-minimal solution $\boldsymbol{x}$ via \hyperlink{target:MinimalS}{$minimal(\boldsymbol{x})$}; hence, steps 22 and 23 remain necessary to ensure that no redundant elements are added to $\mathcal{C}^*$. Excluding $\boldsymbol{y}$ from $\mathcal{X}$ leaves all theoretical results and proofs unchanged.

\subsection{Anytime disjunction algorithm}

The algorithms presented in 3.2 and 3.3 have the property that their output is only meaningful if they complete their tasks in their entirety. Next, we  present an {\it anytime} version.  The algorithm is based on standard disjunctive programming and incorporates one binary variable and two constraints for each generated solution. We present the algorithm with a stopping criterion based on given $\epsilon$ but may be changed and $\epsilon$ may be part of the solution.\\  

Let $\mathcal{C}  \subseteq \mathcal{X}$ and let $\mathcal{Q}^+(\mathcal{C})$ be a disjunctive problem in ($\boldsymbol{x},\boldsymbol{z},\epsilon$) defined as follows:
\begin{flalign*}
	\max&\;\;\epsilon&\;&&\mathcal{Q}^+(\mathcal{C})&\\
	s.t.&\;\;\boldsymbol{a}^+(\boldsymbol{x})^T\boldsymbol{y} \geq \underline{\mathbf{a}}^T \boldsymbol{x}+\epsilon&\mbox{ if }\boldsymbol{z}(\boldsymbol{y})=1\;\forall \boldsymbol{y} \in \mathcal{C}\\ &\;\;\boldsymbol{c}^+(\boldsymbol{x})^T\boldsymbol{y} \geq \underline{\mathbf{c}}^T \boldsymbol{x}+\epsilon&\mbox{ if }\boldsymbol{z}(\boldsymbol{y})=0\;\forall \boldsymbol{y} \in \mathcal{C}&\\
	&\;\;\boldsymbol{x} \in \mathcal{X},\;\boldsymbol{z}(\boldsymbol{y}) \in \lbrace 0,1 \rbrace\;\forall \boldsymbol{y} \in \mathcal{C},\,\epsilon \in \mathbb{N} \cup \lbrace 0 \rbrace&\
\end{flalign*}

Since $\mathcal{C}$ is fixed in problem $\mathcal{Q}^+(\mathcal{C})$ we have that $\mathcal{Q}^+(\mathcal{C})$ is a disjunctive 0-1-MILP problem.\\
 
\begin{remark}\label{remark:Qmas}
	Let $\mathcal{C}  \subseteq \mathcal{X}$ then:
	\begin{enumerate}
		\item If $v(\mathcal{Q}^+(\mathcal{C})) = \epsilon$ then $\mathcal{C}$ is an $\epsilon$-$\mathcal{U}$-coverage. 
		\item If $\lbrace [\boldsymbol{y}]:\boldsymbol{y} \in \mathcal{C} \rbrace \subseteq \mathcal{M}$ and $v(\mathcal{Q}^+(\mathcal{C})) = 0$ then $\mathcal{C}$ is a minimum-cardinality $\mathcal{U}$-coverage.
	\end{enumerate} 
\end{remark}

\textbf{Proof}:
\begin{enumerate}
	\item Let us suppose that $v(\mathcal{Q}^+(\mathcal{C})) = \epsilon$ then $(\boldsymbol{x},\boldsymbol{z},\epsilon+1) \notin F(\mathcal{Q}^+(\mathcal{C}))$ for all $\boldsymbol{x} \in \mathcal{X},\;\forall \boldsymbol{z} \in {\lbrace 0,1 \rbrace}^{|\mathcal{C}|}$. Hence, we have that for all $\boldsymbol{x} \in \mathcal{X}$ there exists $\boldsymbol{y} \in \mathcal{C}$ such that $\boldsymbol{a}^+(\boldsymbol{x})^T\boldsymbol{y} \leq \underline{\mathbf{a}}^T \boldsymbol{x}+\epsilon$ and $\boldsymbol{c}^+(\boldsymbol{x})^T\boldsymbol{y} \leq \underline{\mathbf{c}}^T \boldsymbol{x}+\epsilon$, therefore $\mathcal{C}$ is an $\epsilon$-$\mathcal{U}$-coverage.
	\item Let us suppose that $\lbrace [\boldsymbol{y}]:\boldsymbol{y} \in \mathcal{C} \rbrace \subseteq \mathcal{M}$ and $v(\mathcal{Q}^+(\mathcal{C})) = 0$ then from (1) we have that  $\mathcal{C}$ is a $\mathcal{U}$-coverage. Since $\lbrace [\boldsymbol{y}]:\boldsymbol{y} \in \mathcal{C} \rbrace \subseteq \mathcal{M}$ and $\mathcal{C}$ is a $\mathcal{U}$-coverage we have that $\lbrace [\boldsymbol{y}]:\boldsymbol{y} \in \mathcal{C} \rbrace = \mathcal{M}$ and then  $\mathcal{C}$ is a minimum-cardinality $\mathcal{U}$-coverage $\blacksquare$ 
\end{enumerate}
\begin{algorithm}[H]
	\caption{$Anyt(\epsilon)$}
	\hypertarget{target:anytime}{}
	\renewcommand{\thealgorithm}{}
	\begin{algorithmic}[1]
		\State Let $\boldsymbol{x} \in \mathcal{X}$
		\State Use \hyperlink{target:MinimalS}{$minimal(\boldsymbol{x})$} to obtain $\boldsymbol{y} \in \mathcal{X}$ such that $\boldsymbol{y} \prec \boldsymbol{x}$ and $[\boldsymbol{y}] \in \mathcal{M}$
		\State $\mathcal{C} \leftarrow \lbrace \boldsymbol{y} \rbrace$ 
		\State $stop \leftarrow 0$
		\While{$stop=0$}
		\State  solve $\mathcal{Q}^+(\mathcal{C})$
		\If {$v(\mathcal{Q}^+(\mathcal{C})) \leq \epsilon$} 
		\State $stop \leftarrow 1$
		\Else
		\State Let $\boldsymbol{x}$ be an optimal solution for $\mathcal{Q}^+(\mathcal{C})$,
		\State Use \hyperlink{target:MinimalS}{$minimal(\boldsymbol{x})$} to obtain $\boldsymbol{y} \in \mathcal{X}$ such that $\boldsymbol{y} \prec \boldsymbol{x}$ and $[\boldsymbol{y}] \in \mathcal{M}$
		\State $\mathcal{C} \leftarrow \mathcal{C} \cup \lbrace\boldsymbol{y} \rbrace$ 
		\EndIf
		\EndWhile
		\State \textbf{return} $\mathcal{C}$
	\end{algorithmic}
\end{algorithm}

\begin{remark}\label{remark:algorithmanytime}
	$\;$
	\begin{enumerate}
	\item \hyperlink{target:anytime}{$Anyt(\epsilon)$} is finite.
	\item Let $\mathcal{C}$ be the output of \hyperlink{target:anytime}{$Anyt(\epsilon)$}	then:
	\begin{enumerate}
		\item If $\epsilon=0$ then $\mathcal{C}$ is a minimum-cardinality $\mathcal{U}$-coverage
		\item If $\epsilon \in \mathbb{N} \cup \lbrace 0 \rbrace$ then $\mathcal{C}$ is an $\epsilon$-$\mathcal{U}$-coverage. 
	\end{enumerate}
\end{enumerate}
\end{remark}

\textbf{Proof}: 
\begin{enumerate}
	\item \hyperlink{target:anytime}{$Anyt(\epsilon)$} generates a monotonically non-increasing sequence of values. Every solution $\boldsymbol{y}$ generated at steps 2 and 11 satisfies $[\boldsymbol{y}] \in \mathcal{M}$ by \textbf{Remark \ref{remark:minimal2}}. If $[\boldsymbol{y}] \in \{[\boldsymbol{w}] : \boldsymbol{w} \in \mathcal{C}\}$ before step 12 then $v(\mathcal{Q}^+(\mathcal{C})) = 0$. Since $\mathcal{M}$ is a finite set the algorithm finds $\mathcal{C}$ with $v(\mathcal{Q}^+(\mathcal{C})) \leq \epsilon$ in a finite number of iterations.  
	\item Let $\mathcal{C}$ be the output of \hyperlink{target:anytime}{$Anyt(\epsilon)$} then:
	\begin{enumerate}
		\item If $\epsilon = 0$ then from \textbf{Remark \ref{remark:Qmas}} $\mathcal{C}$ is a minimum cardinality $\mathcal{U}$-coverage.
	    \item If $\epsilon \in \mathbb{N}$ then from \textbf{Remark \ref{remark:Qmas}} $\mathcal{C}$ is an $\epsilon$-$\mathcal{U}$-coverage $\blacksquare$  
	\end{enumerate} 
\end{enumerate}

\subsection{Remarks on the disjunctive problems}

The disjunctive programming problems defined above, which must be solved repeatedly to execute the algorithms, are: $\mathcal{S}(\boldsymbol{x})$, $\mathcal{Q}(\mathcal{C},\epsilon)$ and $\mathcal{Q}^+(\mathcal{C})$. Some solvers (such as CPLEX, for example) include tools for solving disjunctive programming problems. Our experience using {\it indicators}  with CPLEX to handle disjunctive constraints proved unsatisfactory for the problems generated in this work. Consequently, we resorted to standard "Big-M" formulations. In all three cases, defining appropriate "Big-M" values requires bounding the expressions $\underline{\mathbf{a}}^T \boldsymbol{x} - \boldsymbol{a}^+(\boldsymbol{x})^T \boldsymbol{y}$ and $\underline{\mathbf{c}}^T \boldsymbol{x} - \boldsymbol{c}^+(\boldsymbol{x})^T \boldsymbol{y}$ for a given $\boldsymbol{y}$ and $\boldsymbol{x} \in \mathcal{X}$. In \textbf{Remark \ref{remark:bounds}}, we present a sequence of increasingly tight bounds, which require progressively greater computational effort. In our experiments, we simply used $\underline{\mathbf{a}}^T \mathbf{1}$ and $\underline{\mathbf{c}}^T \mathbf{1}$. Experimentation with tighter bounds has been deferred for the time being.

\subsection{Parametric algorithms}

Let $\epsilon \in \mathbb{N} \cup \lbrace 0 \rbrace$, let $\mathcal{C}=\lbrace {\boldsymbol{y}^i \rbrace}_1^L \subseteq \mathcal{X}$, let $A(\epsilon) \in {\lbrace 0,1 \rbrace}^{L \times L}$ defined as follows: $A_{ij}(\epsilon) = 1$ if and only if $\boldsymbol{y}^j \prec_\epsilon \boldsymbol{y}^i$, let $\mathcal{W}(\mathcal{C},\epsilon)$ be a set covering  problem in $\boldsymbol{x}$ defined as follows:  
\begin{flalign*}
	\min&\;\;\sum_{j \in [L]} \boldsymbol{x}_j &\;&&\mathcal{W}(\mathcal{C},\epsilon)&\\
	s.t.&\;\;\sum_{j \in [L]} A_{ij}(\epsilon) \boldsymbol{x}_j \geq 1&\forall i \in [L]\\
	&\;\;\boldsymbol{x} \in {\lbrace 0,1 \rbrace}^L,&\
\end{flalign*}

and let $\mathcal{C}(\boldsymbol{x}) = \lbrace \boldsymbol{y}^j:\boldsymbol{x}_j=1,j\in[L]\rbrace$ for all $\boldsymbol{x} \in F(\mathcal{W}(\mathcal{C},\epsilon))$. 

\begin{remark}\label{remark:parametric}
	Let $\epsilon \in \mathbb{N} \cup \lbrace 0 \rbrace$ and let $\mathcal{C}^*$ be a minimum-cardinality $\mathcal{U}$-coverage, let  $\boldsymbol{x}^*$ be an optimal solution for $\mathcal{W}(\mathcal{C}^*,\epsilon)$ then: $\mathcal{C}(\boldsymbol{x}^*)$ is a minimum-cardinality $\epsilon$-$\mathcal{U}$-coverage.
\end{remark}
\textbf{Proof}:
\begin{enumerate}
	\item If $\boldsymbol{x} \in F(\mathcal{W}(\mathcal{C}^*,\epsilon))$ then for all $i \in [L]$ there exists $j \in [L]$ such that $A_{ij}(\epsilon)\boldsymbol{x}_j=1$ and then for all $i \in [L]$ there exists $\boldsymbol{y}^j \in \mathcal{C}(\boldsymbol{x})$ such that $\boldsymbol{y}^j \prec_\epsilon \boldsymbol{y}^i$, therefore  $\mathcal{C}(\boldsymbol{x})$ is an $\epsilon$-$\mathcal{U}$-coverage for all $\boldsymbol{x} \in F(\mathcal{W}(\mathcal{C}^*,\epsilon))$. 
	
	\item Let $\mathcal{H}(\epsilon)$ be an $\epsilon$-$\mathcal{U}$-coverage. Since $\mathcal{C}^*$ is a $\mathcal{U}$-coverage then for all $\boldsymbol{z} \in \mathcal{H}(\epsilon)$ there exists $j(\boldsymbol{z}) \in [L]$ such that $\boldsymbol{y}^{j(\boldsymbol{z})}  \prec \boldsymbol{z}$. Hence,  $\lbrace \boldsymbol{y}^{j(\boldsymbol{z})}:\boldsymbol{z} \in \mathcal{H}(\epsilon)\rbrace$ is an $\epsilon$-$\mathcal{U}$-coverage. Let $\boldsymbol{x} \in {\lbrace 0,1 \rbrace}^L$ with $\boldsymbol{x}_j=1$ if and only if $j \in \lbrace j(\boldsymbol{z}):\boldsymbol{z} \in \mathcal{H}(\epsilon)\rbrace$, we have that $\boldsymbol{x} \in F(\mathcal{W}(\mathcal{C}^*,\epsilon))$ and   $\mathcal{C}(\boldsymbol{x}) = \lbrace \boldsymbol{y}^{j(\boldsymbol{z})}:\boldsymbol{z} \in \mathcal{H}(\epsilon)\rbrace$ with $|\mathcal{C}(\boldsymbol{x})| \leq |\mathcal{H}(\epsilon)|$.
	\item From (1) and (2) we have that if $\boldsymbol{x}^*$ is an optimal solution for $\mathcal{W}(\mathcal{C}^*,\epsilon)$ then $\mathcal{C}(\boldsymbol{x}^*)$ is a minimum-cardinality  $\epsilon$-$\mathcal{U}$-coverage $\blacksquare$
\end{enumerate}

\begin{remark}\label{remark:approxparametric}
	Let $\epsilon \in \mathbb{N} \cup \lbrace 0 \rbrace$ and let $\mathcal{C}(\epsilon^*)$ be an $\epsilon^*$-$\mathcal{U}$-coverage, let  $\boldsymbol{x}^*$ be an optimal solution for $\mathcal{W}(\mathcal{C}(\epsilon^*),\epsilon)$ then: $\mathcal{C}(\boldsymbol{x}^*)$ is an  $(\epsilon+\epsilon^*)$-$\mathcal{U}$-coverage.
\end{remark}

\textbf{Proof}:
If $\boldsymbol{x} \in F(\mathcal{W}(\mathcal{C}(\epsilon^*),\epsilon))$ then for all $i \in [L]$ there exists $j \in [L]$ such that $A_{ij}(\epsilon)\boldsymbol{x}_j=1$ and then for all $i \in [L]$ there exists $\boldsymbol{y}^j \in \mathcal{C}(\boldsymbol{x})$ such that $\boldsymbol{y}^j \prec_\epsilon \boldsymbol{y}^i$, therefore  $\mathcal{C}(\boldsymbol{x})$ is an $(\epsilon^*+\epsilon)$-$\mathcal{U}$-coverage for all $\boldsymbol{x} \in F(\mathcal{W}(\mathcal{C}(\epsilon^*),\epsilon))\;\blacksquare$.\\

Below we present a parametric algorithm. 

\begin{algorithm}[H]
	\caption{$Parametric(\mathcal{C})$}
	\hypertarget{target:Parametric}{}
	\renewcommand{\thealgorithm}{}
	\begin{algorithmic}[1]
		\State Let $\widehat{\boldsymbol{y}} \in \mathcal{C}$. Let  $\widehat{\epsilon} = \min \lbrace \epsilon:\widehat{\boldsymbol{y}} \prec_\epsilon \boldsymbol{y}\; \forall \boldsymbol{y} \in \mathcal{C} \rbrace$. Let \textbf{solution}$(\epsilon) = \lbrace \widehat{\boldsymbol{y}} \rbrace$  and let $k(\epsilon) = 1$ for all $\epsilon \in \mathbb{N} \cup \lbrace 0 \rbrace$ such that $\epsilon \geq \widehat{\epsilon}$ 
		\While {$\widehat{\epsilon} \geq 1$}  
		\State Let $\delta = \widehat{\epsilon}-1$
		\State Compute $A_{ij}(\delta)$
		\State Solve $\mathcal{W}(\mathcal{C},\delta)$ and let $\boldsymbol{x}$ be an optimal solution.
		\State Let  $\widehat{\epsilon} = \min \lbrace  \epsilon: \boldsymbol{x} \in F(\mathcal{W}(\mathcal{C},\epsilon)) \rbrace$
		\State Let $\mathcal{C}(\boldsymbol{x}) = \lbrace \boldsymbol{y}^j:\boldsymbol{x}_j=1,j\in[L]\rbrace$	
		\State Let \textbf{solution}$(\epsilon) = \mathcal{C}(\boldsymbol{x})$ for all $\epsilon \in \mathbb{N} \cup \lbrace 0 \rbrace$ such that $\widehat{\epsilon} \leq \epsilon \leq \delta$
		\State Let $k(\epsilon)=v(\mathcal{W}(\mathcal{C},\delta))$ for all $\epsilon \in \mathbb{N} \cup \lbrace 0 \rbrace$ such that $\widehat{\epsilon} \leq \epsilon \leq \delta$
		\EndWhile
		\State \textbf{return} $k(\epsilon)$ and \textbf{solution}$(\epsilon)$ for all $\epsilon \in \mathbb{N} \cup \lbrace 0 \rbrace$
	\end{algorithmic}
\end{algorithm}

Let $\epsilon \in \mathbb{N} \cup \lbrace 0 \rbrace$, let $\mathcal{H}^*(\epsilon)$ be a minimum-cardinality $\epsilon$-$\mathcal{U}$-coverage and let $k^*(\epsilon) = |\mathcal{H}^*(\epsilon)|$. Note that $k^*(\epsilon)$ is a step function with $k^*(0) = |\mathcal{C}^*|$ and $k^*(\infty) = 1$.

\begin{remark}\label{remark:parametricalg}
$\;$
	\begin{enumerate}
\item Let $\mathcal{C}^*$ be a minimum-cardinality $\mathcal{U}$-coverage and let \textbf{solution}$(\epsilon)$ and $k(\epsilon)$  for all $\epsilon \in \mathbb{N} \cup \lbrace 0 \rbrace$ be the output of \hyperlink{target:Parametric}{$Parametric(\mathcal{C}^*)$} then \textbf{solution}$(\epsilon)$ is a minimum-cardinality $\epsilon$-$\mathcal{U}$-coverage and $k(\epsilon)=k^*(\epsilon)$ for all $\epsilon \in \mathbb{N} \cup \lbrace 0 \rbrace$. 
\item Let $\mathcal{C}(\epsilon^*)$ be an $\epsilon^*$-$\mathcal{U}$-coverage and let \textbf{solution}$(\epsilon)$ and $k(\epsilon)$  for all $\epsilon \in \mathbb{N} \cup \lbrace 0 \rbrace$ be the output of \hyperlink{target:Parametric}{$Parametric(\mathcal{C}(\epsilon^*))$} then \textbf{solution}$(\epsilon)$ is an  $(\epsilon+\epsilon^*)$-$\mathcal{U}$-coverage and $k(\epsilon) \geq k^*(\epsilon^*+\epsilon)$ for all $\epsilon \in \mathbb{N} \cup \lbrace 0 \rbrace$. 
\end{enumerate}
\end{remark}
\textbf{Proof}:
\begin{enumerate}
 \item  Let $\widehat{\boldsymbol{y}} \in \mathcal{C}^*$ and let $\widehat{\epsilon} = \min \lbrace \epsilon:\widehat{\boldsymbol{y}} \prec_\epsilon \boldsymbol{y}\; \forall \boldsymbol{y} \in \mathcal{C}^* \rbrace$. Since  $\mathcal{C}^*$ is a $\mathcal{U}$-coverage, it follows that \textbf{solution}$(\epsilon) = \lbrace \widehat{\boldsymbol{y}} \rbrace$ is a minimum-cardinality  $\epsilon$-$\mathcal{U}$-coverage for all $\epsilon \in \mathbb{N} \cup \lbrace 0 \rbrace$ such that $\epsilon \geq \widehat{\epsilon}$. Let $\delta = \widehat{\epsilon}-1$ and let $\boldsymbol{x}$ be an optimal solution for  $\mathcal{W}(\mathcal{C}^*, \delta)$. From \textbf{Remark \ref{remark:parametric}} we have that \textbf{solution}$(\delta) = \mathcal{C}(\boldsymbol{x})$ is a minimum-cardinality $\delta$-$\mathcal{U}$-coverage. Let  $\widehat{\epsilon} = \min \lbrace  \epsilon: \boldsymbol{x} \in F(\mathcal{W}(\mathcal{C}^*,\epsilon)) \rbrace$. Note that $\boldsymbol{x}$ remains optimal for all $\epsilon \in \mathbb{N} \cup \lbrace 0 \rbrace$ such that $\widehat{\epsilon} \leq \epsilon \leq \delta$ and then  \textbf{solution}$(\delta) = \mathcal{C}(\boldsymbol{x})$ is a minimum-cardinality $\epsilon$-$\mathcal{U}$-coverage for all $\epsilon \in \mathbb{N} \cup \lbrace 0 \rbrace$ such that $\widehat{\epsilon} \leq \epsilon \leq \delta$.  Next we set $\delta=\widehat{\epsilon}-1$ and so on until $\widehat{\epsilon} = 0$.\\
 
 Since \textbf{solution}$(\epsilon)$ is a minimum-cardinality $\epsilon$-$\mathcal{U}$-coverage for all $\epsilon \in \mathbb{N} \cup \lbrace 0 \rbrace$ we have by construction that $k(\epsilon) = k^*(\epsilon)$ for all $\epsilon \in \mathbb{N} \cup \lbrace 0 \rbrace$.\\
   
\item Let $\widehat{\boldsymbol{y}} \in \mathcal{C}(\epsilon^*)$ and let $\widehat{\epsilon} = \min \lbrace \epsilon:\widehat{\boldsymbol{y}} \prec_\epsilon \boldsymbol{y}\; \forall \boldsymbol{y} \in \mathcal{C}(\epsilon^*) \rbrace$. Since  $\mathcal{C}(\epsilon^*)$ is an $\epsilon^*$-$\mathcal{U}$-coverage, it follows that \textbf{solution}$(\epsilon) = \lbrace \widehat{\boldsymbol{y}} \rbrace$ is an $(\epsilon+\epsilon^*)$-$\mathcal{U}$-coverage for all $\epsilon \in \mathbb{N} \cup \lbrace 0 \rbrace$ such that $\epsilon \geq \widehat{\epsilon}$. Let $\delta = \widehat{\epsilon}-1$ and let $\boldsymbol{x}$ be an optimal solution for  $\mathcal{W}(\mathcal{C}(\epsilon^*), \delta)$. From \textbf{Remark \ref{remark:approxparametric}} we have that \textbf{solution}$(\delta) = \mathcal{C}(\boldsymbol{x})$ is an $(\delta+\epsilon^*)$-$\mathcal{U}$-coverage. Let  $\widehat{\epsilon} = \min \lbrace  \epsilon: \boldsymbol{x} \in F(\mathcal{W}(\mathcal{C}(\epsilon^*),\epsilon)) \rbrace$. Note that $\boldsymbol{x}$ remains optimal for all $\epsilon \in \mathbb{N} \cup \lbrace 0 \rbrace$ such that $\widehat{\epsilon} \leq \epsilon \leq \delta$ and then  \textbf{solution}$(\delta) = \mathcal{C}(\boldsymbol{x})$ is an $(\epsilon+\epsilon^*)$-$\mathcal{U}$-coverage for all $\epsilon \in \mathbb{N} \cup \lbrace 0 \rbrace$ such that $\widehat{\epsilon} \leq \epsilon \leq \delta$. Next we set $\delta=\widehat{\epsilon}-1$ and so on until $\widehat{\epsilon} = 0$.\\

Since \textbf{solution}$(\epsilon)$ is an $(\epsilon^*+\epsilon)$-$\mathcal{U}$-coverage for all $\epsilon \in \mathbb{N} \cup \lbrace 0 \rbrace$ we have by construction that $k(\epsilon) \geq k^*(\epsilon^*
+\epsilon)$ for all $\epsilon \in \mathbb{N} \cup \lbrace 0 \rbrace\;\blacksquare$\\
\end{enumerate}

It should be noted that several solutions with the same cardinality may define the same step of $k^*(\epsilon)$  ($k(\epsilon)$). To simplify the definition of $k^*(\epsilon)$ ($k(\epsilon)$), it suffices to present, for each step, the solution that is feasible at the leftmost point.\\

If $\mathcal{C}(\epsilon^*)$ was generated using \hyperlink{target:anytime}{$Anyt(\epsilon^*)$}, then there exists a sequence of sets $\{\mathcal{C}^l\}_1^K$ such that $v(\mathcal{Q}^+(\mathcal{C}^{l+1})) \leq v(\mathcal{Q}^+(\mathcal{C}^l))$ for all $l \in [K-1]$, and $v(\mathcal{Q}^+(\mathcal{C}^K)) \leq \epsilon^*$. Note that $\mathcal{C}^l$ is a $v(\mathcal{Q}^+(\mathcal{C}^l))$-$\mathcal{U}$-coverage. In general, as $\epsilon^*+\epsilon$ approaches $\epsilon^*$ within the parameterization, it may occur that $\epsilon^*+\epsilon > v(\mathcal{Q}^+(\mathcal{C}^l))$ for coverages of the same cardinality; in such cases, we use $\mathcal{C}^l$ instead of the coverage generated by \hyperlink{target:Parametric}{$Parametric(\mathcal{C}(\epsilon^*))$}.

\section{Computational experience}

The experiments  have been performed on a personal computer as follows: Intel(R)Core(TM) i5-7200U CPU, @ 2.50 GHz Lenovo Zodiac ThinkPad T570, 32.00 GB Ram and Windows 10 Pro Operating System. All the instances have been processed through ILOG-Cplex 22.1.0 from a DOcplex Python code. All the parameters of ILOG-Cplex 22.1.0 are in their default values, except for the $mip.tolerance.integrality$ for which we use $10^{-12}$, for the  required $mip.tolerances.gap$  for \hyperlink{target:anytime}{$Anyt(\epsilon)$} for which we use 0\% and for   $emphasis.mip$ for \hyperlink{target:Disjunctive}{$Disj(\epsilon)$} and \hyperlink{target:Branch and Driven-Disjunction}{$Driven(\epsilon)$} for which we use  1.\\

The computational experiment is conducted using three types of bi-objective problems: the set covering ($\mathcal{SC}$) problem, the selection ($\mathcal{SEL}$) problem and the shortest path ($\mathcal{SP}$) problem. A data generator is defined for each problem type. These generators use 5 parameters for $\mathcal{SC}$ and $\mathcal{SEL}$ problems, and 4 parameters for $\mathcal{SP}$ problems. The data generators are used to control the dimensions of the problems as well as some of the fundamental characteristics of the data.\\

Fourteen experiments related to the construction of a minimum-cardina-
lity $\mathcal{U}$-coverage  were defined. Experiments 1 through 5 (one for each parameter) address $\mathcal{SC}$-type problems; experiments 6 through 10 (one for each parameter) deal with $\mathcal{SEL}$-type problems; and experiments 11 through 14 (one for each parameter) concern $\mathcal{SP}$-type problems. In each experiment, denoted as $e$ in the tables, we fix the parameters used to generate the data and subsequently vary one of them. For each value of the parameter being modified, we generate 10 random instances of the problem. Each set of 10 instances is denoted as $s$ in the tables. To identify a specific set in our subsequent analysis, we simply refer to the pair $(e,s)$.\\

Some sets were used as a reference for varying parameters and appear in multiple experiments to facilitate the comparison of results; in such cases, the set appears several times in the tables, associated with different experiments, and, from its second appearance onwards, we denote it as $(s)$ instead of $s$. For example in Tables \ref{SCparameters} and \ref{SCavemax}, set 3 ($s=3$) appears in experiment 1 ($e=1$) and reappears (as $s=(3)$) in experiments 3 ($e=3$) and 5 ($e=5$).\\

A total of 16 batches were generated for $\mathcal{SC}$-type problems, 14 batches for $\mathcal{SEL}$-type problems, and 10 batches for $\mathcal{SP}$-type problems, resulting in a total of 40 batches and 400 instances. For each of the 400 generated instances, we determined a $\mathcal{U}$-coverage using \hyperlink{target:Disjunctive}{$Disj(0)$}, \hyperlink{target:Branch and Driven-Disjunction}{$Driven(0)$}, or \hyperlink{target:anytime}{$Anyt(0)$}.\\

For each of the 400 generated instances, a parametric construction of a minimum-cardinality $\epsilon$-$\mathcal{U}$-coverage is performed, for all $\epsilon \in \mathbb{N} \cup \{0\}$, based on the previously constructed minimum-cardinality $\mathcal{U}$-coverage.\\

Let $\epsilon^* = \lfloor 0.01 \times \min\lbrace \underset{\boldsymbol{x} \in \mathcal{X}}\min\; \underline{\mathbf{a}}^T \boldsymbol{x},\underset{\boldsymbol{x} \in \mathcal{X}}\min\; \underline{\mathbf{c}}^T \boldsymbol{x} \rbrace \rfloor$. We use $\epsilon^*$ to perform an approximate parametric analysis. One batch is selected from each experiment for each problem type; for every instance within that batch a parametric construction of a $(\epsilon^*+\epsilon)$-$\mathcal{U}$-coverage is carried out, based on an $\epsilon^*$-$\mathcal{U}$-coverage found by using \hyperlink{target:Disjunctive}{$Disj(\epsilon^*)$}, \hyperlink{target:Branch and Driven-Disjunction}{$Driven(\epsilon^*)$} or \hyperlink{target:anytime}{$Anyt(\epsilon^*)$}.\\

In addition to the parameters for generating the instances, the tables present, for each $(e,s)$ and, where applicable, for each specifically identified algorithm, the average and maximum values of the following statistics:
\begin{itemize}
	\item  $|\mathcal{C}^*|$: with $\mathcal{C}^*$ a minimum-cardinality $\mathcal{U}$-coverage,  
	\item $|PF|$: with $PF$ the Pareto frontier for   $(\underline{\mathbf{a}},\underline{\mathbf{c}})$,
	\item \#$\mathcal{T}$: number of $\mathcal{T}$-type problems solved to find a minimum-cardinality $\mathcal{U}$-coverage, 
	\item \#$\mathcal{RS}$: number of $\mathcal{R}$-type and $\mathcal{S}$-type problems solved to find a minimum-cardinality $\mathcal{U}$-coverage,
	\item $t$: time, in seconds, to find a minimum-cardinality $\mathcal{U}$-coverage,
	\item $t_p$: time, in seconds, to perform a parametric analysis,
	\item $t_g$: time, in seconds, to find an $\epsilon^*$-$\mathcal{U}$-coverage and to perform an approximate parametric analysis,
	\item $ps$: number of steps to define $k^*(\epsilon)$, 
	\item $k*$: minimum number of solutions used to define an $\epsilon$-$\mathcal{U}$-coverage, with $* \in \lbrace 15,10,5,2,1 \rbrace$ and
	\[\epsilon \leq \lfloor\frac{*}{100} \times  \min\lbrace \underset{\boldsymbol{x} \in \mathcal{X}}\min\; \underline{\mathbf{a}}^T \boldsymbol{x},\underset{\boldsymbol{x} \in \mathcal{X}}\min\; \underline{\mathbf{c}}^T \boldsymbol{x} \rbrace \rfloor\] according to the method presented in the table.  
\end{itemize} 

It should be noted that the number of problems of type $\mathcal{Q}$ or type $\mathcal{Q}^+$ solved to find a minimum-cardinality $\mathcal{U}$-coverage  is equal to $|\mathcal{C}^*|$.\\

For the three types of problems,  in the corresponding initial experiment (experiment 1 for $\mathcal{SC}$, 6 for $\mathcal{SEL}$, and 11 for $\mathcal{SP}$), we varied the problem dimensions ($n$ or $|V|$, as applicable) and evaluated the performance of the three algorithms (\hyperlink{target:Branch and Driven-Disjunction}{$Driven(0)$},  \hyperlink{target:Disjunctive}{$Disj(0)$} and  \hyperlink{target:anytime}{$Anyt(0)$}). Based on the information obtained from these experiments, we selected one of the three algorithms to evaluate while varying the other parameters.\\

The results will be presented separately for each type of problem, with the exception of those related to parameterization. For each case, the problem, the data generator, and the experiments concerning the construction of a minimum-cardinality $\mathcal{U}$-coverage will be described. The results for exact and approximate parameterizations are analyzed together, as the observations are consistent.\\

In Tables \ref{SCave}, \ref{SCmax}, \ref{SELave}, and \ref{SELmax}, a "$-$" appears in the columns corresponding to $\#\mathcal{RS}$ and $t$ for certain pairs $(e,s)$. This indicates that, for the set $s$, the experiment was not conducted using the specified algorithm. This decision was made because the growth trend of the processing time for these algorithms became very pronounced as the problem dimensions increased. Nevertheless, all these cases were processed by at least one of the algorithms.

\subsection{Set covering problems}

\subsubsection{Problem and data generator} 

Let $n \geq 1$. Let $\mathbf{A} \in \lbrace 0,1 \rbrace^{n \times n}$. The  bi-objective set covering problem parameterized by  $(\boldsymbol{a},\boldsymbol{c})$ considered, denoted as $\mathcal{SC}(\boldsymbol{a},\boldsymbol{c})$, is defined as follows: 	  
\begin{flalign*}
	\min&\;\;\sum_{j \in  [n]}\boldsymbol{a}_j\boldsymbol{x}_j,\sum_{j \in  [n]}\boldsymbol{c}_j\boldsymbol{x}_j&\;&&\mathcal{SC}(\boldsymbol{a},\boldsymbol{c})\\
	s.t.&\;\;\sum_{j \in [n]} \mathbf{A}_{ij}\boldsymbol{x}_j\geq 1&\forall i \in [n]&\\
	&\;\;\boldsymbol{x} \in \lbrace 0,1 \rbrace^n&\\
\end{flalign*}

Let $n \geq 1$, let $\beta,\gamma,\delta \in (0,1]$ and let $\overline{d} \geq 1000$.\\ 

In order to generate $\mathbf{A}$ we proceed as follows:

\begin{itemize}
	\item let $\mathbf{A}_{ij} = 0$ for all $(i,j) \in [n] \times [n]$,
	\item generate $n$ random points in $[0,1] \times [0,1]$ and let $d_{ij}$ the distance from point $i$ to point $j$,
	\item let $maxd = \max\lbrace d_{ij}:i,j \in[n] \rbrace$
	\item if $d_{ij} < \beta \times maxd$ we make $\mathbf{A}_{ij} = 1$ for all $i,j \in [n]$.
\end{itemize}

In order to generate uncorrelated  data we proceed as follows:
\begin{itemize}
	\item generate $\underline{\mathbf{a}}_j$ by taking a random integer in $[1000,\overline{d}]$ for all $j \in [n]$,
	\item let $\overline{\mathbf{a}} = \underline{\mathbf{a}}$,
	\item with probability $\delta$ let $\overline{\mathbf{a}}_j = \lfloor (1+u\gamma)\underline{\mathbf{a}}_j \rfloor$ for all $j \in [n]$ with $u$ taken randomly in $(0,1)$,
	\item proceed analogously to generate $\underline{\mathbf{c}}$ and $\overline{\mathbf{c}}$.
\end{itemize}

Table \ref{SCparameters} shows the parameters used to generate the instances. There are 16 instance sets, numbered 1 through 16. For each $(e,s)$, where $e = 1, \dots, 5$ and $s = 1, \dots, 16$, the table presents the values $(n, \beta, \gamma, \delta, \overline{d})$, as well as the average and maximum cardinality of the Pareto frontier associated with $(\underline{\mathbf{a}}, \underline{\mathbf{c}})$.

\subsubsection{Experiments 1 through 5 regarding the construction of a  minimum-cardinality $\mathcal{U}$-coverage}

In Experiment 1, we set $(\beta, \gamma, \delta, \overline{d}) = (0.15, 0.05, 1.0, 20000)$ and vary $n$ across the set $\{100, 150, 200, 250, 300, 350\}$ to define sets 1 through 6. Table \ref{SCave} clearly shows that: (i) the average execution times for the algorithms increase as $n$ grows, as would be expected; (ii) the average times for \hyperlink{target:Branch and Driven-Disjunction}{$Driven(0)$} and   \hyperlink{target:Disjunctive}{$Disj(0)$} are significantly lower than those for \hyperlink{target:anytime}{$Anyt(0)$}; and (iii) the average times for \hyperlink{target:Branch and Driven-Disjunction}{$Driven(0)$} from $n=250$ are lower than those for \hyperlink{target:Disjunctive}{$Disj(0)$}. Table \ref{SCmax} clearly shows that: (i) the worst-case execution times for the algorithms in each set increase as $n$ grows, as would be expected; (ii) the worst-case times for \hyperlink{target:Branch and Driven-Disjunction}{$Driven(0)$} and   \hyperlink{target:Disjunctive}{$Disj(0)$} are significantly lower than those for \hyperlink{target:anytime}{$Anyt(0)$} ; and (iii) the worst-case times for \hyperlink{target:Branch and Driven-Disjunction}{$Driven(0)$} from $n=250$ are significantly lower than those for \hyperlink{target:Disjunctive}{$Disj(0)$}.\\

The point is that, while it is true, as one might expect, that $\#\mathcal{T} > |\mathcal{C}^*|$ and that the value of $\#\mathcal{RS}$ associated with \hyperlink{target:Branch and Driven-Disjunction}{$Driven(0)$} is greater than those associated with \hyperlink{target:Disjunctive}{$Disj(0)$} and \hyperlink{target:anytime}{$Anyt(0)$}, the significantly lower computational effort required for problems of type $\mathcal{T}$, compared to disjunctive problems of type $\mathcal{Q}$ and $\mathcal{Q}^+$ solved using \hyperlink{target:Disjunctive}{$Disj(0)$} and \hyperlink{target:anytime}{$Anyt(0)$}, clearly tips the balance in favor of \hyperlink{target:Branch and Driven-Disjunction}{$Driven(0)$}.\\

Based on the results of Experiment 1, we selected \hyperlink{target:Branch and Driven-Disjunction}{$Driven(0)$} to continue the experiments for $\mathcal{SC}$.\\

In Experiment 2, we use set 4 as a reference labeled as (4). The parameter $\beta$ is adjusted, by both decreasing and increasing it,  to define sets 7 and 8. The parameter $\beta$ influences the density of the coverage matrix $\mathbf{A}$. The results presented in Table \ref{SCavemax} indicate, as expected, a strong dependence of the performance of \hyperlink{target:Branch and Driven-Disjunction}{$Driven(0)$} with respect to the density of the matrix.\\

In Experiment 3, we use set 3 as a reference, labeled as (3). The parameter $\gamma$ is increased to define sets 9 and 10. This parameter influences the deviation of the data. The results presented in Table \ref{SCavemax} indicate, as expected, a strong dependence of the  performance of \hyperlink{target:Branch and Driven-Disjunction}{$Driven(0)$} with respect to the increase in the uncertainty range.\\

In Experiment 4, we use set 4 as a reference, labeled as (4). The parameter $\delta$ is reduced to define sets 11, 12, and 13. The parameter $\delta$ represents the probability that each element of the objectives is incremented. The results presented in Table \ref{SCavemax} indicate, as expected, a strong dependence of the  performance of \hyperlink{target:Branch and Driven-Disjunction}{$Driven(0)$} with respect to the decrease in the probability of movement within the objective data.\\

In Experiment 5, we use set 3 as a reference, labeled with (3). The parameter $\overline{d}$ is reduced to define sets 14, 15, and 16. Since $\overline{d}$ acts as an upper bound for defining $\underline{\mathbf{a}}$ and $\underline{\mathbf{c}}$, its reduction causes the nominal data values to become closer to one another. The results presented in Table \ref{SCavemax} indicate, as might be expected, a strong dependence of the performance of  \hyperlink{target:Branch and Driven-Disjunction}{$Driven(0)$} with respect to the decrease in $\overline{d}$; however, in this case, the performance does not exhibit monotonicity.

\subsection{Selection problem}

\subsubsection{Problem and data generator}

Let $n,p \geq 1$ with $p \leq n$. The bi-objective $n,p$-selection problem   parameterized by  $(\boldsymbol{a},\boldsymbol{c})$ considered, denoted as $\mathcal{SEL}(\boldsymbol{a},\boldsymbol{c})$, is defined as follows:
\begin{flalign*}
	\min&\;\;\sum_{j \in [n]} \boldsymbol{a}_j\boldsymbol{x}_j,\sum_{j \in [n]}  \boldsymbol{c}_j\boldsymbol{x}_j&\;&&\mathcal{SEL}(\boldsymbol{a},\boldsymbol{c})&\\
	s.t.&\;\;\sum_{j \in [n]} \boldsymbol{x}_j = p&\\
	&\;\;\boldsymbol{x} \in \lbrace 0,1 \rbrace^n&\\
\end{flalign*}

Let $n,p \geq 1$ with $p \leq n$ and let $\gamma,\delta \in (0,1]$ and let $\overline{d} \geq 1000$.\\ 

Table \ref{SELparameters} shows the parameters used to generate the instances. There are 14 problem sets, numbered 1 through 14. For each $(e,s)$ where $e = 6, \dots, 10$ and $s = 1, \dots, 14$ the table presents the values $(n, p, \gamma, \delta, \overline{d})$, as well as the average and maximum cardinality of the Pareto front associated with $(\underline{\mathbf{a}}, \underline{\mathbf{c}})$.

\subsubsection{Experiments 6 through 10 regarding the construction of a  minimum-cardinality $\mathcal{U}$-coverage}

In Experiment 6, we set $(\gamma, \delta, \overline{d}) = (0.05, 1.0, 20000)$ and vary $(n, p)$ across the set $\{(50, 5), (80, 8), (100, 10), (120, 12)\}$ to define sets 1 through 4. Tables \ref{SELave} and \ref{SELmax} clearly show that: (i) the average and worst-case execution times for the algorithms increase, as expected, as $(n, p)$ grows; (ii) the average and worst-case times for \hyperlink{target:Branch and Driven-Disjunction}{$Driven(0)$} and   \hyperlink{target:Disjunctive}{$Disj(0)$} are significantly lower than those for \hyperlink{target:anytime}{$Anyt(0)$}; and (iii) the average and worst-case times for \hyperlink{target:Branch and Driven-Disjunction}{$Driven(0)$} are significantly lower than those for \hyperlink{target:Disjunctive}{$Disj(0)$}\\

What happens, once again, is that, while it is to be expected that $\#\mathcal{T} > |\mathcal{C}^*|$ and $\#\mathcal{RS}$ associated with \hyperlink{target:Branch and Driven-Disjunction}{$Driven(0)$} is greater than that associated with \hyperlink{target:Disjunctive}{$Disj(0)$} and \hyperlink{target:anytime}{$Anyt(0)$},  the significantly lower computational effort required for problems of type $\mathcal{T}$, compared to disjunctive problems of type $\mathcal{Q}$ and $\mathcal{Q}^+$ solved using \hyperlink{target:Disjunctive}{$Disj(0)$} and \hyperlink{target:anytime}{$Anyt(0)$}, clearly tips the balance in favor of \hyperlink{target:Branch and Driven-Disjunction}{$Driven(0)$}.\\

Based on the results of Experiment 6, we selected \hyperlink{target:Branch and Driven-Disjunction}{$Driven(0)$} to continue the experiments for $\mathcal{SEL}$.\\

In Experiment 7, we use set 3 as reference, labeled with (3). The parameter $p$ is adjusted (increased or decreased) to define sets 5 and 6. The results presented in Table \ref{SELavemax} indicate, as one might expect, a strong dependence of the \hyperlink{target:Branch and Driven-Disjunction}{$Driven(0)$} performance on the value of $p$.\\

In Experiment 8, we once again use set 3 as reference, labeled with (3). The parameter $\gamma$ is increased to define sets 7 and 8; this parameter influences the deviation of the data. The results presented in Table \ref{SELavemax} indicate, as expected, a strong dependence of the \hyperlink{target:Branch and Driven-Disjunction}{$Driven(0)$} performance on the increase in the range of uncertainty.\\

In Experiment 9, we use set 4 as reference, labeled with (4). The parameter $\delta$ is reduced to define sets 9, 10, and 11. The parameter $\delta$ represents the probability that each element of the objectives increases. The results presented in Table \ref{SELavemax} indicate a dependency, but not a monotonicity, in the performance of \hyperlink{target:Branch and Driven-Disjunction}{$Driven(0)$} as the probability of movement in the objective data decreases.\\

In Experiment 10, we use set 3 as reference, labeled with (3). The parameter $\overline{d}$ is reduced to define  sets 12, 13, and 14. Since $\overline{d}$ acts as an upper bound for the definition of $\underline{\mathbf{a}}$ and $\underline{\mathbf{c}}$, its reduction causes the nominal data values to become closer to one another. The results presented in Table \ref{SELavemax} show no significant changes in \hyperlink{target:Branch and Driven-Disjunction}{$Driven(0)$} performance as $\overline{d}$ decreases; it may be worth considering more drastic modifications to identify any underlying trends.

\subsection{Shortest path problem}

\subsubsection{Problem and data generator}

Consider a directed graph without \emph{loops} $G=(V,E)$ with $|V| \geq 1$. Let $\hat v \in V$ and $\hat w \in V$. We are looking for paths from $\hat v$ to $\hat w$. For each $v \in V$ let $E^+(v) = \lbrace (v,w): (v,w) \in E  \rbrace$ and let $E^-(v) = \lbrace (w,v):(w,v) \in E \rbrace$. Let $\boldsymbol{x} \in {\lbrace 0,1 \rbrace}^{|E|}$ with $\boldsymbol{x}_e = 1$ if the edge $e$ is selected to define a path. The bi-objective shortest path  problem parameterized by $(\boldsymbol{a},\boldsymbol{c})$ considered, denoted as $\mathcal{SP}(\boldsymbol{a},\boldsymbol{c})$, is defined as follows:
\begin{flalign*}
	\min&\;\;\sum_{e \in E} \boldsymbol{a}_e\boldsymbol{x}_e,\sum_{e \in E} \boldsymbol{c}_e\boldsymbol{x}_e&\;&&\mathcal{SP}(\boldsymbol{a},\boldsymbol{c})&\\
	s.t.&\;\;\sum_{e:e \in E^+(\hat v)} \boldsymbol{x}_e = 1,\;\;\;\sum_{e:e \in E^-(\hat w)} \boldsymbol{x}_e = 1&\\
	&\;\;\sum_{e:e \in E^-(v)} \boldsymbol{x}_e  - \sum_{e \in E^+(v)} \boldsymbol{x}_e =0\;\;\forall v \in V\;\; \mbox{with}\;\; v \notin \lbrace \hat v,\hat w \rbrace\\
	&\;\;\mathbf{x} \in {\lbrace 0,1 \rbrace}^{|E|}\\
\end{flalign*}

The graphs  are  generated following \cite{Hanasusanto2015} as follows:
\begin{itemize}
	\item in each problem instance, the nodes correspond to $|V|$ points with coordinates  that are chosen uniformly at random in the square $[0,1] \times [0,1]$,
	\item we choose the pair of nodes with the largest
	Euclidean distance as $\hat v$ and $\hat w$ (the start and terminal nodes),
	\item in order to generate $E$ we begin with a fully connected graph and remove 70\% of the arcs  in order of decreasing Euclidean distance,  that is, starting with the longest arcs.
\end{itemize}

Let $|V| \geq 1$, let $\gamma,\delta \in (0,1]$ and let $\overline{d} \geq 1000$.\\

In order to generate correlated data we proceed as follows:
\begin{itemize}
	\item generate $\mathbf{P}_i$ by taking a random integer in $[1000,\overline{d}]$ for all $i \in [|V|]$.
	\item let $\underline{\mathbf{a}}_e = \lfloor 10000 \times d_{ij}\rfloor$ with $e$ the edge from node $i$ to node $j$ and $d_{ij}$ the Euclidean distance from node $i$ to node $j$.
	\item with probability $\delta$ we make $\overline{\mathbf{a}}_e = \lfloor (1+u\gamma)\underline{\mathbf{a}}_e \rfloor$ for all $e \in E$ with $u$ taken randomly in $(0,1)$.
	\item We make $\underline{\mathbf{c}}_e =  |\mathbf{P}_i - \mathbf{P}_j|$  with $e$ the edge  from node $i$ to node $j$
	\item With probability $\delta$ we make $\overline{\mathbf{c}}_e = \lfloor (1+u\gamma)\underline{\mathbf{c}}_e \rfloor$ for all $e \in E$ with $u$ taken randomly in $(0,1)$.
\end{itemize}

Since all components $\underline{\mathbf{a}}_e$ are positive in the generated data we have that all generated minimal elements correspond to simple paths and then it is not necessary to include constraints to avoid cycles.\\

Table \ref{SPparameters} shows the parameters used to generate the instances. There are 10 problem sets, numbered 1 through 10. For each $(e,s)$ where $e = 11, \dots, 14$ and $s = 1, \dots, 10$ the table presents the values $(|V|,n=|E|,\gamma, \delta, \overline{d})$, as well as the average and maximum cardinality of the Pareto front associated with $(\underline{\mathbf{a}}, \underline{\mathbf{c}})$.

\subsubsection{Experiments 11 through 14 regarding the construction of a  minimum-cardinality $\mathcal{U}$-coverage}

In Experiment 11, we set $(\gamma, \delta, \overline{d}) = (0.05, 1.0, 20000)$ and vary $|V|$ across the set $\{30, 40, 50\}$ to define sets 1 through 3. As can be clearly observed in Tables \ref{SPave} and \ref{SPmax}: (i) the average and worst-case running times of the algorithms increase as $|V|$ grows, which is to be expected; (ii) the average and worst-case running times for \hyperlink{target:Branch and Driven-Disjunction}{$Driven(0)$} and \hyperlink{target:Disjunctive}{$Disj(0)$}  are clearly lower than those for  \hyperlink{target:anytime}{$Anyt(0)$}; and (iii) somewhat surprisingly, the average and worst-case running times for \hyperlink{target:Disjunctive}{$Disj(0)$} are clearly lower than those for \hyperlink{target:Branch and Driven-Disjunction}{$Driven(0)$}.\\

What is happening in this case is that, given that $\#\mathcal{T}$ is much larger than $|\mathcal{C}^*|$, and the $\#\mathcal{RS}$ associated with \hyperlink{target:Branch and Driven-Disjunction}{$Driven(0)$} is much larger than the $\#\mathcal{RS}$ associated with \hyperlink{target:Disjunctive}{$Disj(0)$} and \hyperlink{target:anytime}{$Anyt(0)$}, the lower computational effort required for $\mathcal{T}$-type problems compared to the disjunctive $\mathcal{Q}$-type and $\mathcal{Q}^+$-type problems solved using \hyperlink{target:Disjunctive}{$Disj(0)$} and \hyperlink{target:anytime}{$Anyt(0)$} does not sufficiently offset the overall effort, ultimately favoring the use of \hyperlink{target:Disjunctive}{$Disj(0)$}.\\

Admittedly, one could argue that the choice of the second objective is rather peculiar (even though it has been used to evaluate algorithms designed for shortest-path problems); nevertheless, this case illustrates that dismissing disjunctive programming for multi-objective problems is not always the appropriate approach.\\

Based on the results of Experiment 11 we selected \hyperlink{target:Disjunctive}{$Disj(0)$} to continue the experiments for $\mathcal{SP}$.\\

In Experiment 12, we use set 2 as reference, labeled (2). The parameter $\gamma$ is increased to define sets 4 and 5. This parameter influences the deviation of the data.  As expected, the results shown in Table \ref{SPavemax} indicate a strong dependence of the \hyperlink{target:Disjunctive}{$Disj(0)$} performance on the increase in the uncertainty range.\\

In Experiment 13, we use 3 as reference, labeled (3). The parameter $\delta$ is decreased to define sets 6 and 7. The parameter $\delta$ represents the probability that each element of the objectives may or may not increase. The results shown in Table \ref{SPavemax} indicate a dependency, but not monotonicity, in the performance of  \hyperlink{target:Disjunctive}{$Disj(0)$} as the probability of movement in the objective data decreases.\\

In Experiment 14, we use set 2 as reference, labeled (2). The parameter $\overline{d}$ is decreased to define sets 8, 9, and 10. The parameter $\overline{d}$ serves as the upper bound for defining $\underline{\mathbf{a}}$ and $\underline{\mathbf{c}}$; as it decreases, the nominal values of the data move closer together. The results shown in Table \ref{SPavemax} show an abrupt change when the parameter reaches the value 5000; prior to this, the performance was very consistent. It may be worth considering more  modifications to identify any underlying trends.

\subsection{Parametric and approx parametric analysis}

Tables \ref{SCparametric}, \ref{SELparametric}, and \ref{SPparametric} show the results of the exact parametric analysis, based on prior knowledge of a minimum-cardinality $\mathcal{U}$-coverage, for the 400 generated instances. For each pair $(e,s)$, we present the average values and the corresponding worst-case values for the execution time of the parameterization algorithm, the number of segments defining the function $k^*(\epsilon)$, as well as the values $k15,k10,k5,k2,k1$.\\

An approximate parametric analysis was performed on a selected set from each experiment using two algorithms: the one used with $\epsilon = 0$, now running with $\epsilon^*$ and \hyperlink{target:anytime}{$Anyt(\epsilon^*)$}. This can be seen in the tables \ref{SCapproxparametric}, \ref{SELapproxparametric}, and \ref{SPapproxparametric}. For reference and to facilitate comparisons, the results with $\epsilon = 0$ are presented. Averages and worst-case scenarios are shown for each $(e,s)$, for total run time (time to find coverage plus parameterization time), and for the values of $k15,k10,k5,k2,k1$.\\ 

The times associated with the parameterizations are clearly correlated with the times to find the $\mathcal{U}$-coverages or the $\epsilon^*$-$\mathcal{U}$-coverages. Observing the columns corresponding to the exact values of $k15,k10,k5,k2,k1$, it is notable that, for the generated instances, it seems possible to maintain an inventory of solutions with which to calculate the Pareto frontier of the revealed scenario. Likewise, the approximate parameterizations provide, in much less time, an inventory of solutions with which the approximations for $k15,k10,k5,k2,k1$ are more than acceptable. It is worth noting that approximations based on \hyperlink{target:anytime}{$Anyt(\epsilon^*)$} have better quality as the cardinality of the coverages approaches the exact value of $k1$.

\section{Conclusions}

In this paper, we use an extension of the concepts of coverage and $\epsilon$-coverage from deterministic bi-objective combinatorial optimization to bi-objective combinatorial optimization problems with integer data under box uncertainty ($\mathcal{U}$-$BOCO$). By employing extreme scenario characterizations on the uncertainty box, we define preorders on the feasible set whose induced posets allow us to identify a minimum-cardinality $\mathcal{U}$-coverage and an $\epsilon$-$\mathcal{U}$-coverage by finding minimal elements.\\ 

We propose three exact algorithmic frameworks to construct $\epsilon$-$\mathcal{U}$-covera-
ges and minimum-cardinality $\mathcal{U}$-coverages: a disjunctive programming algorithm (\hyperlink{target:Disjunctive}{$Disj(\epsilon)$}) that adds one binary variable and two constraints per generated solution, a driven-disjunctive algorithm (\hyperlink{target:Branch and Driven-Disjunction}{$Driven(\epsilon)$}) that avoids introducing additional binary variables through an implicit enumeration scheme over disjunctions, and an anytime disjunctive algorithm (\hyperlink{target:anytime}{$Anyt$}) that may treat $\epsilon$ as a decision variable (even though, in the computational experiments, we have used a pre-set $\epsilon$ to define the stopping criterion for the algorithm).\\

We present exact and approximate parametric algorithms based on set-covering formulations to characterize the trade-off between coverage quality $\epsilon$ and  coverage cardinality.\\

Computational experiments across 400 instances of set covering, selection, and shortest path problems show that the proposed methods are effective for moderate problem sizes. Notably, no single exact algorithm dominates across all problem structures. Furthermore, the parametric analysis shows that while \hyperlink{target:anytime}{$Anyt(0)$} is generally slower for exact full-coverage generation, the anytime approach and approximate parameterization may remain valuable when time limits apply or when an $\epsilon$-approximate coverage is acceptable. It is important to note that for the same predetermined $\epsilon$, the quality of the parametric analysis based on \hyperlink{target:anytime}{$Anyt(\epsilon)$} was superior, on the generated data, to that based on the other algorithms.\\

Future work includes investigating tighter Big-$M$ bounds for the disjunctive subproblems and extending the methodology to problems with more than two objectives, further strengthening the practical scope of order-theoretic and disjunctive approaches under uncertainty.

\clearpage
\begin{table}[h]
	\caption{$\mathcal{SC}$ problems. Parameters for data generation.}\label{SCparameters}
	\begin{tabular*}{\textwidth}{@{\extracolsep\fill}ccccccccc}
\noalign{\smallskip}\hline\noalign{\smallskip}			
$e$&$s$&$n$&$\beta$&$\gamma$&$\delta$&$\overline{d}$&$ave|PF|$&$max|PF|$\\
\midrule
1&	1&100&0.150&0.050&1.0&20000&37.9&56\\
&	2&150&&&&&70.0&115\\
&	3&200&&&&&67.0&92\\
&	4&250&&&&&85.5&137\\
&	5&300&&&&&90.3&156\\
&	6&350&&&&&93.1& 162\\
\noalign{\smallskip}\hline\noalign{\smallskip}	
2&	7&250&0.125&0.050&1.0&20000&131.7&181\\
&	(4)&  &0.150&&&&85.5&137\\
&	8&  &0.175&&&&66.9&105\\
\noalign{\smallskip}\hline\noalign{\smallskip}
3&	(3)&200&0.150&0.050&1.0&20000&67.0&92\\
&	9&&&0.075\\
&	10&&&0.085\\
\noalign{\smallskip}\hline\noalign{\smallskip}	
4&(4)&250&0.150&0.050&1.0&20000&85.5&137\\
&11&&&&0.9\\
&12&&&&0.8\\
&13&&&&0.6\\
\noalign{\smallskip}\hline\noalign{\smallskip}	
5&(3)&200&0.150&0.050&1.0&20000&67.0&92\\
&14&&&&&10000&73.8&141\\
&15&&&&&7500&70.5&109\\
&16&&&&&5000&47.1&67\\
	\end{tabular*}
\end{table}

\begin{table}[h]
	\caption{$\mathcal{SC}$ problems. Average values.}\label{SCave}
	\begin{tabular*}{\textwidth}{@{\extracolsep\fill}cccccccccc}
		\noalign{\smallskip}\hline\noalign{\smallskip}
		&  \multicolumn{2}{@{}c@{}}{} & \multicolumn{3}{@{}c@{}}{\hyperlink{target:Branch and Driven-Disjunction}{$Driven(0)$}} & \multicolumn{2}{@{}c@{}}{\hyperlink{target:Disjunctive}{$Disj(0)$}}  & \multicolumn{2}{@{}c@{}}{\hyperlink{target:anytime}{$Anyt(0)$}}
		\\\cmidrule{4-6}\cmidrule{7-8}\cmidrule{9-10}%
		$e$&$s$&$|\mathcal{C}^*|$&\#$\mathcal{T}$&\#$\mathcal{RS}$&$t$&\#$\mathcal{RS}$&$t$&\#$\mathcal{RS}$&$t$\\
		\noalign{\smallskip}\hline\noalign{\smallskip}
		1&1&62.8&140.2&111.4&29.3&96.9&14.7&62.8&33.4\\
		&2&133.3& 332.4& 255.8& 148.6&196.5&136.6&133.3&305.6\\
		&3&124.8& 322.4& 256.3 &228.0&187.7&186.2&124.8&454.9\\
		&4&199.8&583.0& 466.0& 691.4&291.1&1137.9&199.8&2622.3\\
		&5&174.5&430.0&337.9&530.5&258.6&884.8&$-$&$-$\\
		&6&196.8& 557.4& 445.7& 958.5&$-$&$-$&$-$&$-$\\
	\end{tabular*}      
\end{table}

\begin{table}[h]
	\caption{$\mathcal{SC}$ problems. Maximum values.}\label{SCmax}
	\begin{tabular*}{\textwidth}{@{\extracolsep\fill}cccccccccc}
		\noalign{\smallskip}\hline\noalign{\smallskip}
		&  \multicolumn{2}{@{}c@{}}{} & \multicolumn{3}{@{}c@{}}{\hyperlink{target:Branch and Driven-Disjunction}{$Driven(0)$}} & \multicolumn{2}{@{}c@{}}{\hyperlink{target:Disjunctive}{$Disj(0)$}}  & \multicolumn{2}{@{}c@{}}{\hyperlink{target:anytime}{$Anyt(0)$}}		
		\\\cmidrule{4-6}\cmidrule{7-8}\cmidrule{9-10}%
		$e$&$s$&$|\mathcal{C}^*|$&\#$\mathcal{T}$&\#$\mathcal{RS}$&$t$&\#$\mathcal{RS}$&$t$&\#$\mathcal{RS}$&$t$\\
		\noalign{\smallskip}\hline\noalign{\smallskip}
		1&1&111&256&195&66&168&34&111&110\\
		&2&203& 566& 447& 284&283&286&203&817\\
		&3&232& 656& 496& 542&319&589&232&1616\\
		&4&381 &1396& 1136& 2213&506&4005&381&9869\\
		&5&368&960&725&1648&512&5477&$\;\;-$&$\;\;-$\\
		&6&353& 1112& 937& 2402&$\;\;-$&$\;\;-$&$\;\;-$&$\;\;-$\\
	\end{tabular*}      
\end{table}

\begin{table}[h]
	\caption{$\mathcal{SC}$ problems. \protect\hyperlink{target:Branch and Driven-Disjunction}{$Driven(0)$} }\label{SCavemax}
	\begin{tabular*}{\textwidth}{@{\extracolsep\fill}cccccccccc}
		\noalign{\smallskip}\hline\noalign{\smallskip}
		&  \multicolumn{1}{@{}c@{}}{} & \multicolumn{4}{@{}c@{}}{ave} & \multicolumn{4}{@{}c@{}}{max}
		\\\cmidrule{3-6}\cmidrule{7-10}%
		$e$&$s$&$|\mathcal{C}^*|$&\#$\mathcal{T}$&\#$\mathcal{RS}$&$t$&$|\mathcal{C}^*|$&\#$\mathcal{T}$&\#$\mathcal{RS}$&$t$\\
		\noalign{\smallskip}\hline\noalign{\smallskip}
		2&7&379.5&1184.8&922.5&2125.4&728&2500&1949&6701\\
		&(4)&199.8&583.0& 466.0& 691.4&381 &1396& 1136& 2213\\
		&8&118.9&275.0&215.0&226.3&206&456&354&426\\
		\noalign{\smallskip}\hline\noalign{\smallskip}		
		3&(3)&124.8& 322.4& 256.3 &228.0&232& 656& 496& 542\\
		&9&181.4&607.4&488.1&532.4&345&1782&1516&2037\\
		&10&213.1&756.6&815.0&604.2&430&2192&1843&3056\\
		\noalign{\smallskip}\hline\noalign{\smallskip}			
		4&(4)&199.8&583.0& 466.0& 691.4 &381&1396& 1136& 2213\\
		&11&186.9&508.4&403.5&591.9&343&1064&839&1588\\
		&12&168.6&415.2&332.0&442.9&314&808&620&1071\\
		&13&143.6&328.6&262.4&306.5&291&710&536&882\\
		\noalign{\smallskip}\hline\noalign{\smallskip}
		5&(3)&124.8& 322.4& 256.3 &228.0&232& 656& 496& 542\\	
		&14&187.8&547.4&424.6&462.8&429&1302&972&1553\\
		&15&194.7&694.6&562.8&640.1&360&2044&1841&2286\\
		&16&129.0&371.4&284.2&237.5&204&654&505&504\\	
	\end{tabular*}      
\end{table}

\begin{table}[h]
	\caption{$\mathcal{SEL}$ problems. Parameters for data generation}\label{SELparameters}
	\begin{tabular*}{\textwidth}{@{\extracolsep\fill}ccccccccc}
		\noalign{\smallskip}\hline\noalign{\smallskip}
		$e$&$s$&$n$&$p$&$\gamma$&$\delta$&$\overline{d}$&$ave|PF|$&$max|PF|$\\
		\noalign{\smallskip}\hline\noalign{\smallskip}
		6&1&50&5&0.050&1.0&20000&30.0&45\\
		&2&80&8&&&&58.5&123	\\
		&3&100&10&&&&93.7&142\\
		&4&120&12&&&&156.1&264\\
		\noalign{\smallskip}\hline\noalign{\smallskip}	
		7&5&100&15&0.050&1.0&20000&159.2&307\\
		&(3)&&10&&&&93.7&142\\	
		&6&&5&&&&39.7&61\\	
		\noalign{\smallskip}\hline\noalign{\smallskip}
		8&(3)&100&10&0.050&1.0&20000&93.7&142\\
		&7&&&0.075&\\
		&8&&&0.100&\\	
		\noalign{\smallskip}\hline\noalign{\smallskip}
		9&(4)&120&12&0.050&1.0&20000&156.1&264\\
		&9&&&&0.9&\\
		&10&&&&0.8&\\	
		&11&&&&0.6&\\
		\noalign{\smallskip}\hline\noalign{\smallskip}	
		10&(3)&100&10&0.050&1.0&20000&93.7&142\\
		&12&&&&&10000&106.3&153\\
		&13&&&&&7500&87.7&141\\
		&14&&&&&5000&102.2&157\\
	\end{tabular*}      
\end{table}

\begin{table}[h]
	\caption{$\mathcal{SEL}$ problems. Average values.}\label{SELave}
	\begin{tabular*}{\textwidth}{@{\extracolsep\fill}cccccccccc}
		\noalign{\smallskip}\hline\noalign{\smallskip}
				&  \multicolumn{2}{@{}c@{}}{} & \multicolumn{3}{@{}c@{}}{\hyperlink{target:Branch and Driven-Disjunction}{$Driven(0)$}} & \multicolumn{2}{@{}c@{}}{\hyperlink{target:Disjunctive}{$Disj(0)$}}  & \multicolumn{2}{@{}c@{}}{\hyperlink{target:anytime}{$Anyt(0)$}}
		\\\cmidrule{4-6}\cmidrule{7-8}\cmidrule{9-10}%
		$e$&$s$&$|\mathcal{C}^*|$&\#$\mathcal{T}$&\#$\mathcal{RS}$&$t$&\#$\mathcal{RS}$&$t$&\#$\mathcal{RS}$&$t$\\
		\noalign{\smallskip}\hline\noalign{\smallskip}
		6&1&37.4&78.4&67.5&4.1&61.6&2.5&37.4&5.3\\
		&2&92.1&201.6&162.9&29.0&146.7&24.9&92.1&76.4\\	
		&3&163.5&399.0&318.5&93.2&251.4&137.3&163.5&339.5\\	
		&4&340.8&856.4&652.8&431.3&490.7&4002.2&$-$&$-$\\
	\end{tabular*}      
\end{table}

\begin{table}[h]
	\caption{$\mathcal{SEL}$ problems. Maximum values.}\label{SELmax}
	\begin{tabular*}{\textwidth}{@{\extracolsep\fill}cccccccccc}
		\noalign{\smallskip}\hline\noalign{\smallskip}
		&  \multicolumn{2}{@{}c@{}}{} & \multicolumn{3}{@{}c@{}}{\hyperlink{target:Branch and Driven-Disjunction}{$Driven(0)$}} & \multicolumn{2}{@{}c@{}}{\hyperlink{target:Disjunctive}{$Disj(0)$}}  & \multicolumn{2}{@{}c@{}}{\hyperlink{target:anytime}{$Anyt(0)$}}
		\\\cmidrule{4-6}\cmidrule{7-8}\cmidrule{9-10}%
		$e$&$s$&$|\mathcal{C}^*|$&\#$\mathcal{T}$&\#$\mathcal{RS}$&$t$&\#$\mathcal{RS}$&$t$&\#$\mathcal{RS}$&$t$\\
		\noalign{\smallskip}\hline\noalign{\smallskip}
		6&1&61&138&111&9&90&6&61&18\\	
		&2&188&404&315&78&293&83&188&261\\
		&3&301&1008&823&343&434&506&301&1484\\
		&4&787&2424&1852&1895&1061&26517&$-$&$-$\\
	\end{tabular*}      
\end{table}

\begin{table}[h]
	\caption{$\mathcal{SEL}$ problems. \protect\hyperlink{target:Branch and Driven-Disjunction}{$Driven(0)$} }\label{SELavemax}
	\begin{tabular*}{\textwidth}{@{\extracolsep\fill}cccccccccc}
		\noalign{\smallskip}\hline\noalign{\smallskip}
		&  \multicolumn{1}{@{}c@{}}{} & \multicolumn{4}{@{}c@{}}{ave} & \multicolumn{4}{@{}c@{}}{max}
		\\\cmidrule{3-6}\cmidrule{7-10}%
		$e$&$s$&$|\mathcal{C}^*|$&\#$\mathcal{T}$&\#$\mathcal{RS}$&$t$&$|\mathcal{C}^*|$&\#$\mathcal{T}$&\#$\mathcal{RS}$&$t$\\
		\noalign{\smallskip}\hline\noalign{\smallskip}
		7&5&370.8&985.6&757.2&570.0&948&3132&2426&3273\\	
		&(3)&163.5&399.0&318.5&93.2&301&1008&823&343\\	
		&6&57.0&121.6&102.9&14.4&107&232&180&37\\
		\noalign{\smallskip}\hline\noalign{\smallskip}
		8&(3)&163.5&399.0&318.5&93.2&301&1008&823&343\\
		&7&219.3&537.6&394.6&152.5&470&1190&813&545\\
		&8&273.5&750.4&554.9&264.5&615&1886&1367&1083\\
		\noalign{\smallskip}\hline\noalign{\smallskip}	
		9&(4)&340.8&856.4&652.8&431.3&787&2424&1852&1895\\
		&9&339.9&1327.2&1140.1&893.5&1054&7826&7091&7083\\
		&10&319.0&871.6&700.9&415.0&815&3314&2778&2377\\	
		&11&258.5&603.6&493.4&214.0&480&1272&1066&610\\	
		\noalign{\smallskip}\hline\noalign{\smallskip}
		10&(3)&163.5&399.0&318.5&93.2&301&1008&823&343\\
		&12&225.6&551.6&422.2&141.8&392&958&702&349\\
		&13&184.7&474.6&376.7&113.8&326&954&809&269\\
		&14&275.7&768.0&584.7&240.7&514&1550&1221&634\\	
	\end{tabular*}      
\end{table}

\begin{table}[h]
	\caption{$\mathcal{SP}$ problems. Parameters for data generation.}\label{SPparameters}
	\begin{tabular*}{\textwidth}{@{\extracolsep\fill}ccccccccc}
		\noalign{\smallskip}\hline\noalign{\smallskip}
		$e$&$s$&$|V|$&$n=|E|$&$\gamma$&$\delta$&$\overline{d}$&$ave|PF|$&$max|PF|$\\
		\noalign{\smallskip}\hline\noalign{\smallskip}
		11&1&30&262&0.050&1.0&20000&4.2&7\\
		&2&40&468&&&&4.4&7\\
		&3&50&736&&&&6.7&10\\
		\noalign{\smallskip}\hline\noalign{\smallskip}
		12&(2)&40&468&0.050&1.0&20000&4.4&7\\	
		&4&&&0.075&\\
		&5&&&0.100&\\
		\noalign{\smallskip}\hline\noalign{\smallskip}	
		13&(3)&50&736&0.050&1.0&20000&6.7&10\\
		&6&&&&0.9&\\
		&7&&&&0.8&\\
		\noalign{\smallskip}\hline\noalign{\smallskip}	
		14&(2)&40&468&0.050&1.0&20000&4.4&7\\
		&8&&&&&10000&6.5&14\\
		&9&&&&&7500&6.1&13\\
		&10&&&&&5000&3.7&8\\
	\end{tabular*}      
\end{table}

\begin{table}[h]
	\caption{$\mathcal{SP}$ problems. Average values.}\label{SPave}
	\begin{tabular*}{\textwidth}{@{\extracolsep\fill}cccccccccc}
		\noalign{\smallskip}\hline\noalign{\smallskip}
				&  \multicolumn{2}{@{}c@{}}{} & \multicolumn{3}{@{}c@{}}{\hyperlink{target:Branch and Driven-Disjunction}{$Driven(0)$}} & \multicolumn{2}{@{}c@{}}{\hyperlink{target:Disjunctive}{$Disj(0)$}}  & \multicolumn{2}{@{}c@{}}{\hyperlink{target:anytime}{$Anyt(0)$}}
		\\\cmidrule{4-6}\cmidrule{7-8}\cmidrule{9-10}%
		$e$&$s$&$|\mathcal{C}^*|$&\#$\mathcal{T}$&\#$\mathcal{RS}$&$t$&\#$\mathcal{RS}$&$t$&\#$\mathcal{RS}$&$t$\\
		\noalign{\smallskip}\hline\noalign{\smallskip}
		11&1&26.7&72.2&48.0&23.8&29.4&5.1&26.8&8.9\\
		&2&45.5&200.0&157.8&135.5&49.7&20.9&45.5&27.4\\	
		&3&193.7&1165.4*&1081.8*&2536.8*&204.9&719.3&194.7&2004.9\\
	\end{tabular*}
	\footnotetext{*: averages taken over the 7 simplest problems of the set according to the \hyperlink{target:Disjunctive}{$Disj(0)$} execution.}      
\end{table}

\begin{table}[h]
	\caption{$\mathcal{SP}$ problems. Maximum values.}\label{SPmax}
	\begin{tabular*}{\textwidth}{@{\extracolsep\fill}cccccccccc}
		\noalign{\smallskip}\hline\noalign{\smallskip}
				&  \multicolumn{2}{@{}c@{}}{} & \multicolumn{3}{@{}c@{}}{\hyperlink{target:Branch and Driven-Disjunction}{$Driven(0)$}} & \multicolumn{2}{@{}c@{}}{\hyperlink{target:Disjunctive}{$Disj(0)$}}  & \multicolumn{2}{@{}c@{}}{\hyperlink{target:anytime}{$Anyt(0)$}}
		\\\cmidrule{4-6}\cmidrule{7-8}\cmidrule{9-10}%
		$e$&$s$&$|\mathcal{C}^*|$&\#$\mathcal{T}$&\#$\mathcal{RS}$&$t$&\#$\mathcal{RS}$&$t$&\#$\mathcal{RS}$&$t$\\
		\noalign{\smallskip}\hline\noalign{\smallskip}
		11&1&99&294&197&145&109&27&99&57\\
		&2&116&680&566&599&126&99&116&143\\
		&3&583&4080*&3872*&11290*&608&3716&583&13840\\
	\end{tabular*}
	\footnotetext{*: maximums taken over the 7 simplest problems of the set according to the \hyperlink{target:Disjunctive}{$Disj(0)$} execution.}       
\end{table}

\begin{table}[h]
	\caption{$\mathcal{SP}$ problems. \protect\hyperlink{target:Disjunctive}{$Disj(0)$}}\label{SPavemax}
	\begin{tabular*}{\textwidth}{@{\extracolsep\fill}cccccccc}
		\cmidrule{1-8}
		&  \multicolumn{1}{@{}c@{}}{} & \multicolumn{4}{@{}c@{}}{ave} & \multicolumn{2}{@{}c@{}}{max}
		\\\cmidrule{3-6}\cmidrule{7-8}%
		$e$&$s$&$|\mathcal{C}^*|$&\#$\mathcal{RS}$&$t$&$|\mathcal{C}^*|$&\#$\mathcal{RS}$&$t$\\
		\cmidrule{1-8}
		12&(2)&45.5&49.7&20.9&116&126&99\\	
		&4&65.7&70.9&44.9&214&229&262\\	
		&5&89.1&95.4&87.6&290&309&542\\
			\cmidrule{1-8}
		13&(3)&193.7&204.9&719.3&583&608&3716\\
		&6&117.6&131.4&224.6&343&363&1167\\
		&7&71.4&82.0&77.8&188&222&288\\
			\cmidrule{1-8}
		14&(2)&45.5&49.7&20.9&116&126&99\\
		&8&63.3&68.2&38.8&174&186&189\\
		&9&56.8&63.0&26.1&98&111&51\\
		&10&131.1&138.1&1066.6&858&888&10420\\
	\end{tabular*}      
\end{table}

\begin{table}[h]
	\caption{$\mathcal{SC}$ problems. Exact parametric analysis by using \protect\hyperlink{target:Parametric}{$Parametric(\mathcal{C}^*)$}} \label{SCparametric}
	\begin{tabular*}{\textwidth}{@{\extracolsep\fill}cccccccccccccccc}
		\cmidrule{1-16}
		&  \multicolumn{1}{@{}c@{}}{} & \multicolumn{7}{@{}c@{}}{$ave$} & \multicolumn{7}{@{}c@{}}{$max$}
		\\\cmidrule{3-9}\cmidrule{10-16}%
		$e$&$s$&$t_p$&$ps$&$k15$&$k10$&$k5$&$k2$&$k1$&$t_p$&$ps$&$k15$&$k10$&$k5$&$k2$&$k1$\\
		\cmidrule{1-16}
		1&1&5.0&46.1&4.2& 6.3 &11.4 &23.1 &35.4 &16& 85& 8 &11 &22 &41& 64\\
		&2&38.4& 90.0& 5.8& 9.0& 17.9& 41.1& 69.5& 86& 128& 8 &13& 26 &69& 104\\
		&3&35.6& 86.6&  7.0& 10.4&  20.6&  43.6& 68.6&  116&  136&  11&  16&  30&  64&  104\\
		&4&127.3 & 137.0 &9.3& 14.1& 28.4 &68.6 &116.1 &475 &259 &13 &18& 40& 112& 213\\
		&5&83.1&117.2&8.9&13.3&25.6&58.1&96.3&424&204&16&24&53&136&221\\
		&6&107.7& 131.9& 9.3 &14.2& 27.2 &62.5 &102.4& 351 &201 &13 &19 &39& 105& 189\\
		\cmidrule{1-16}
		2&7&658.1& 210.9 &8.4 &12.7& 28.2 &78.1& 154.5 &2599 &318 &13 &20 &47 &144 &286\\
		&(4)&127.3 & 137.0 &9.3& 14.1& 28.4 &68.6 &116.1 &475 &259 &13 &18& 40& 112& 213\\
		&8&24.7&81.3&10.8&15.5&27.9&53.2&79.6&94&126&13&19&36&71&147\\
		\cmidrule{1-16}
		3&		(3)&35.6& 86.6&  7.0& 10.4&  20.6&  43.6& 68.6&  116&  136&  11&  16&  30&  64&  104\\
		&9&78.6&114.1&7.5&11.2&24.0&54.2&95.5&260&198&11&17&36&95&173\\
		&10&129.5&129.0&7.5&11.5&25.0&61.6&109.1&485&247&11&18&39&108&199\\	
		\cmidrule{1-16}
		4&(4)&127.3 & 137.0 &9.3& 14.1& 28.4 &68.6 &116.1 &475 &259 &13 &18& 40& 112& 213\\
		&11&119.5&133.2&9.3&14.1&28.1&65.2&108.1&351&239&13&19&39&97&196\\
		&12&95.1&123.1&9.2&13.6&27.5&59.7&98.5&284&229&13&19&38&90&187\\
		&13&62.8&108.5&9.0&13.5&26.2&50.9&84.2&264&217&13&18&38&79&157\\
		\cmidrule{1-16}
		5&(3)&35.6& 86.6&  7.0& 10.4&  20.6&  43.6& 68.6&  116&  136&  11&  16&  30&  64&  104\\
		&14&96.7&116.1&5.2&8.2&15.9&40.8&80.0&488&206&7&12&21&63&179\\
		&15&85.6&114.0&4.6&7.0&14.2&37.3&71.6&338&189&6&9&19&65&135\\
		&16&21.8&89.1&3.1&4.5&9.6&25.4&48.5&60&125&4&6&12&36&68\\				
	\end{tabular*}      
\end{table}

\begin{table}[h]
	\caption{$\mathcal{SC}$ problems. Approximate parametric analysis by using \protect\hyperlink{target:Parametric}{$Parametric(\mathcal{C}(\epsilon^*))$}. $\epsilon^* = \lfloor 0.01 \times \min\lbrace \underset{\boldsymbol{x} \in \mathcal{X}}\min\; \underline{\mathbf{a}}^T \boldsymbol{x},\underset{\boldsymbol{x} \in \mathcal{X}}\min\; \underline{\mathbf{c}}^T \boldsymbol{x} \rbrace \rfloor$. $\mathcal{C}(\epsilon^*)$ generated with the algorithm indicated}\label{SCapproxparametric}
	\begin{tabular*}{\textwidth}{@{\extracolsep\fill}cccccccccccccccc}
		\cmidrule{1-14}
		&  \multicolumn{1}{@{}c@{}}{} & \multicolumn{6}{@{}c@{}}{$ave$} & \multicolumn{6}{@{}c@{}}{$max$}
		\\\cmidrule{3-8}\cmidrule{9-14}%
		$e(s)$&$alg$&$t_g$&$k15$&$k10$&$k5$&$k2$&$k1$&$t_g$&$k15$&$k10$&$k5$&$k2$&$k1$\\
		\cmidrule{1-14}
		1(6)&\hyperlink{target:Branch and Driven-Disjunction}{$Driven(0)$}&1066.2&  9.3 &14.2& 27.2 &62.5 &102.4&2661 &13 &19 &39& 105& 189\\
		&\hyperlink{target:Branch and Driven-Disjunction}{$Driven(\epsilon^*)$}&414.1&9.7&15.2&32.2&101.7&130.1&985&14&20&51&189&226\\
		& \hyperlink{target:anytime}{$Anyt(\epsilon^*)$}&734.4&9.7&15.3&32.4&81.8&115.2&2263&14&20&51&138&199\\
		\cmidrule{1-14}
		2(7)&\hyperlink{target:Branch and Driven-Disjunction}{$Driven(0)$}&2783.5&  8.4 &12.7& 28.2 &78.1& 154.5 &9300 &13 &20 &47 &144 &286\\
		&\hyperlink{target:Branch and Driven-Disjunction}{$Driven(\epsilon^*)$}&682.0&8.8&14.7&35.4&153.0&216.8&2470&13&24&57&281&457\\
		& \hyperlink{target:anytime}{$Anyt(\epsilon^*)$}&2071.8&8.8&14.7&35.1&103.5&174.9&7962&13&24&57&191&310\\
		\cmidrule{1-14}
		3(10)&\hyperlink{target:Branch and Driven-Disjunction}{$Driven(0)$}&944.5&7.5&11.5&25.0&61.6&109.1&3541&11&18&39&108&199\\
		&\hyperlink{target:Branch and Driven-Disjunction}{$Driven(\epsilon^*)$}&348.7&8.2&12.9&32.6&107.6&150.3&1255&12&20&48&199&315\\
		& \hyperlink{target:anytime}{$Anyt(\epsilon^*)$}&526.6&8.2&12.9&32.7&79.2&125.1&1740&12&20&48&132&223\\
		\cmidrule{1-14}
		4(4) &\hyperlink{target:Branch and Driven-Disjunction}{$Driven(0)$}&823.3  &9.3& 14.1& 28.4 &68.6 &116.1 &2696&13 &18& 40& 112& 213\\
		&\hyperlink{target:Branch and Driven-Disjunction}{$Driven(\epsilon^*)$}&339.1&9.7&15.9&34.3&115.0&142.4&1043&14&22&52&209&276\\
		& \hyperlink{target:anytime}{$Anyt(\epsilon^*)$}&823.5&9.7&15.8&34.3&84.3&125.3&2696&14&22&52&144&231\\
		\cmidrule{1-14}
		5(15) &\hyperlink{target:Branch and Driven-Disjunction}{$Driven(0)$}&725.7&4.6&7.0&14.2&37.3&71.6&2489&6&9&19&65&135\\
		&\hyperlink{target:Branch and Driven-Disjunction}{$Driven(\epsilon^*)$}&164.7&5.0&7.9&18.1&71.2&105.8&445&6&11&26&132&208\\
		& \hyperlink{target:anytime}{$Anyt(\epsilon^*)$}&220.7&5.2&7.8&17.8&51.5&87.1&650&6&10&26&97&151\\
	\end{tabular*}      
\end{table}

\begin{table}[h]
	\caption{$\mathcal{SEL}$ problems. Exact parametric analysis by using \protect\hyperlink{target:Parametric}{$Parametric(\mathcal{C}^*)$} }
	\label{SELparametric}
	\begin{tabular*}{\textwidth}{@{\extracolsep\fill}cccccccccccccccc}
		\cmidrule{1-16}
		&  \multicolumn{1}{@{}c@{}}{} & \multicolumn{7}{@{}c@{}}{$ave$} & \multicolumn{7}{@{}c@{}}{$max$}
		\\\cmidrule{3-9}\cmidrule{10-16}%
		$e$&$s$&$t_p$&$ps$&$k15$&$k10$&$k5$&$k2$&$k1$&$t_p$&$ps$&$k15$&$k10$&$k5$&$k2$&$k1$\\
		\cmidrule{1-16}
		6&1&0.8&28.5&9.9&12.8&19.4&30.0&32.6&3&39&15&22&27&51&52\\	
		&2&19.1&61.5&9.6&14.1&23.8&45.5&64.4&74&109&18&27&46&107&138\\
		&3&66.0&104.4&9.9&15.3&28.3&58.4&91.6&236&160&13&20&37&80&161\\
		&4&625.8&170.8&11.6&17.8&35.0&82.4&162.9&3269&293&17&24&49&137&317\\
		\cmidrule{1-16}
		7&5&1355.6&186.5&8.7&12.7&26.1&60.7&120.1&6411&330&12&17&35&85&243\\		
		&(3)&66.0&104.4&9.9&15.3&28.3&58.4&91.6&236&160&13&20&37&80&161\\
		&6&4.4&40.8&11.7&15.3&23.2&38.8&46.8&17&59&17&21&31&58&74\\
		\cmidrule{1-16}
		8&(3)&66.0&104.4&9.9&15.3&28.3&58.4&91.6&236&160&13&20&37&80&161\\
		&7&133.7&126.2&10.4&16.1&30.6&73.2&132.6&606&202&14&20&39&112&262\\
		&8&230.8&147.0&11.0&16.9&34.8&92.5&170.3&1091&250&14&21&46&158&338\\
		\cmidrule{1-16}
		9&(4)&625.8&170.8&11.6&17.8&35.0&82.4&162.9&3269&293&17&24&49&137&317\\	
		&9&850.9&170.7&11.7&17.5&35.1&81.5&152.2&6514&352&17&24&49&154&369\\
		&10&557.3&160.3&11.6&17.4&34.1&74.9&133.9&3557&286&17&25&47&113&266\\
		&11&262.7&145.4&11.2&17.2&33.2&70.6&122.2&1019&235&17&26&46&110&196\\
		\cmidrule{1-16}
		10&(3)&66.0&104.4&9.9&15.3&28.3&58.4&91.6&236&160&13&20&37&80&161\\
		&12&134.7&125.0&6.9&10.0&20.5&51.0&98.8&439&182&11&16&34&93&170\\
		&13&87.7&102.2&4.9&7.2&14.9&37.0&68.9&244&153&7&11&23&58&119\\
		&14&177.5&116.0&3.6&5.8&11.4&29.8&64.9&595&150&5&7&15&43&81\\
	\end{tabular*}      
\end{table}

\begin{table}[h]
	\caption{$\mathcal{SEL}$ problems. Approximate parametric analysis by using \protect\hyperlink{target:Parametric}{$Parametric(\mathcal{C}(\epsilon^*))$}. $\epsilon^* = \lfloor 0.01 \times \min\lbrace \underset{\boldsymbol{x} \in \mathcal{X}}\min\; \underline{\mathbf{a}}^T \boldsymbol{x},\underset{\boldsymbol{x} \in \mathcal{X}}\min\; \underline{\mathbf{c}}^T \boldsymbol{x} \rbrace \rfloor$.$\mathcal{C}(\epsilon^*)$ generated with the algorithm indicated}\label{SELapproxparametric}
	\begin{tabular*}{\textwidth}{@{\extracolsep\fill}cccccccccccccccc}
		\cmidrule{1-14}
		&  \multicolumn{1}{@{}c@{}}{} & \multicolumn{6}{@{}c@{}}{$ave$} & \multicolumn{6}{@{}c@{}}{$max$}
		\\\cmidrule{3-8}\cmidrule{9-14}%
		$e(s)$&$alg$&$t_g$&$k15$&$k10$&$k5$&$k2$&$k1$&$t_g$&$k15$&$k10$&$k5$&$k2$&$k1$\\
		\cmidrule{1-14}
		6(4)&\hyperlink{target:Branch and Driven-Disjunction}{$Driven(0)$}&1057.1&11.6&17.8&35.0&82.4&162.9&5164&17&24&49&137&317\\
		&\hyperlink{target:Branch and Driven-Disjunction}{$Driven(\epsilon^*)$}&326.9&12.7&19.5&42.5&161.2&214.6&1265&19&28&62&315&454\\
		&\hyperlink{target:anytime}{$Anyt(\epsilon^*)$}&898.3&12.7&19.5&42.0&112.3&183.6&3384&19&28&62&179&352\\
		\cmidrule{1-14}
		7(5)&\hyperlink{target:Branch and Driven-Disjunction}{$Driven(0)$}&1925.6&8.7&12.7&26.1&60.7&120.1&9684&12&17&35&85&243\\
		&\hyperlink{target:Branch and Driven-Disjunction}{$Driven(\epsilon^*)$}&195.9&9.1&14.2&32.4&116.8&172.8&1024&12&19&44&239&396\\
		&\hyperlink{target:anytime}{$Anyt(\epsilon^*)$}&461.5&9.3&14.3&32.3&88.8&149.3&2410&12&19&44&134&307\\
		\cmidrule{1-14}
		8(8)&\hyperlink{target:Branch and Driven-Disjunction}{$Driven(0)$}&495.3&11.0&16.9&34.8&92.5&170.3&2174&14&21&46&158&338\\
		&\hyperlink{target:Branch and Driven-Disjunction}{$Driven(\epsilon^*)$}&252.4&11.8&18.2&43.9&170.0&206.2&1031&15&23&61&337&463\\
		&\hyperlink{target:anytime}{$Anyt(\epsilon^*)$}&568.8&11.9&18.2&43.9&117.7&181.3&2624&15&23&61&215&364\\
		\cmidrule{1-14}
		9(9)&\hyperlink{target:Branch and Driven-Disjunction}{$Driven(0)$}&1744.4&11.7&17.5&35.1&81.5&152.2&13597&17&24&49&154&369\\
		&\hyperlink{target:Branch and Driven-Disjunction}{$Driven(\epsilon^*)$}&449.6&12.6&19.6&41.8&150.5&207.8&2974&19&28&60&361&611\\
		&\hyperlink{target:anytime}{$Anyt(\epsilon^*)$}&1140.4&12.5&19.7&41.6&110.9&177.8&7679&19&28&60&212&434\\
		\cmidrule{1-14}
		10(10)&\hyperlink{target:Branch and Driven-Disjunction}{$Driven(0)$}&972.3&11.6&17.4&34.1&74.9&133.0&5934&17&25&47&113&266\\
		&\hyperlink{target:Branch and Driven-Disjunction}{$Driven(\epsilon^*)$}&259.2&12.4&19.4&41.0&132.1&182.9&1271&19&28&63&269&422\\
		&\hyperlink{target:anytime}{$Anyt(\epsilon^*)$}&665.4&12.3&19.5&40.9&101.8&157.5&3551&19&29&63&160&297\\
	\end{tabular*}      
\end{table}

\begin{table}[h]
	\caption{$\mathcal{SP}$ problems.Exact parametric analysis by using \protect\hyperlink{target:Parametric}{$Parametric(\mathcal{C}^*)$}}\label{SPparametric}
	\begin{tabular*}{\textwidth}{@{\extracolsep\fill}cccccccccccccccc}
		\cmidrule{1-16}
		&  \multicolumn{1}{@{}c@{}}{} & \multicolumn{7}{@{}c@{}}{$ave$} & \multicolumn{7}{@{}c@{}}{$max$}
		\\\cmidrule{3-9}\cmidrule{10-16}%
		$e$&$s$&$t_p$&$ps$&$k15$&$k10$&$k5$&$k2$&$k1$&$t_p$&$ps$&$k15$&$k10$&$k5$&$k2$&$k1$\\
		\cmidrule{1-16}
		11&1&0.6&16.8&1.3&1.6&2.4&4.6&9.3&5&38&2&3&5&12&25\\
		&2&0.9&26.8&1.6&2.0&3.1&8.8&16.3&8&68&2&3&7&29&47\\
		&3&33.1&45.0&1.3&1.8&2.9&8.8&22.3&211&80&2&4&6&14&47\\
		\cmidrule{1-16}
		12&(2)&0.9&26.8&1.6&2.0&3.1&8.8&16.3&8&68&2&3&7&29&47\\
		&4&4.1&35.3&1.6&2.2&4.7&15.9&27.9&35&106&2&4&14&69&117\\
		&5&9.3&43.6&1.8&2.7&7.4&25.4&40.6&73&131&3&6&30&115&180\\
		\cmidrule{1-16}
		13&(3)&33.1&45.0&1.3&1.8&2.9&8.8&22.3&211&80&2&4&6&14&47\\
		&6&6.3&32.0&1.5&1.9&2.9&6.9&14.9&30&64&3&4&6&13&30\\
		&7&2.4&26.4&1.4&1.9&2.9&6.4&9.7&14&48&2&4&6&12&16\\
		\cmidrule{1-16}
		14&(2)&0.9&26.8&1.6&2.0&3.1&8.8&16.3&8&68&2&3&7&29&47\\
		&8&2.7&31.9&2.3&3.3&5.6&14.0&27.1&18&73&4&5&10&49&108\\
		&9&1.5&29.9&2.6&3.1&5.6&13.0&23.3&3&45&4&5&10&21&42\\
		&10&12.2&20.3&1.8&2.4&4.2&8.2&12.9&113&36&3&6&9&22&26\\
	\end{tabular*}      
\end{table}

\begin{table}[h]
	\caption{$\mathcal{SP}$ problems. Approximate parametric analysis by using \protect\hyperlink{target:Parametric}{$Parametric(\mathcal{C}(\epsilon^*))$}. $\epsilon^* = \lfloor 0.01 \times \min\lbrace \underset{\boldsymbol{x} \in \mathcal{X}}\min\; \underline{\mathbf{a}}^T \boldsymbol{x},\underset{\boldsymbol{x} \in \mathcal{X}}\min\; \underline{\mathbf{c}}^T \boldsymbol{x} \rbrace \rfloor$. $\mathcal{C}(\epsilon^*)$ generated with the algorithm indicated}\label{SPapproxparametric}
	\begin{tabular*}{\textwidth}{@{\extracolsep\fill}cccccccccccccccc}
		\cmidrule{1-14}
		&  \multicolumn{1}{@{}c@{}}{} & \multicolumn{6}{@{}c@{}}{$ave$} & \multicolumn{6}{@{}c@{}}{$max$}
		\\\cmidrule{3-8}\cmidrule{9-14}%
		$e(s)$&$alg$&$t_g$&$k15$&$k10$&$k5$&$k2$&$k1$&$t_g$&$k15$&$k10$&$k5$&$k2$&$k1$\\
		\cmidrule{1-14}
		11(3)&\hyperlink{target:Disjunctive}{$Disj(0)$}&752.4&1.3&1.8&2.9&8.8&22.3&3927&2&4&6&14&47\\
		&\hyperlink{target:Disjunctive}{$Disj(\epsilon^*)$}&14.3&1.5&2.1&3.7&21.0&28.9&48&3&4&7&45&62\\
		&\hyperlink{target:anytime}{$Anyt(\epsilon^*)$}&18.3&1.5&2.0&3.5&12.0&28.2&48&3&5&7&26&54\\
		\cmidrule{1-14}
		12(5)&\hyperlink{target:Disjunctive}{$Disj(0)$}&96.9&1.8&2.7&7.4&25.4&40.6&615&3&6&30&115&180\\
		&\hyperlink{target:Disjunctive}{$Disj(\epsilon^*)$}&32.4&1.8&2.9&9.6&40.7&47.2&235&3&8&43&180&192\\
		&\hyperlink{target:anytime}{$Anyt(\epsilon^*)$}&60.9&1.8&2.9&9.6&38.5&45.8&514&3&8&43&180&189\\
		\cmidrule{1-14}
		13(6)&\hyperlink{target:Disjunctive}{$Disj(0)$}&230.9&1.5&1.9&2.9&6.9&14.9&1197&3&4&6&13&30\\
		&\hyperlink{target:Disjunctive}{$Disj(\epsilon^*)$}&8.0&1.5&1.9&3.6&15.1&19.4&19&3&4&7&30&35\\
		&\hyperlink{target:anytime}{$Anyt(\epsilon^*)$}&11.0&1.5&1.9&3.4&10.3&19.6&19&3&4&7&17&32\\
		\cmidrule{1-14}
		14(10)&\hyperlink{target:Disjunctive}{$Disj(0)$}&1078.8&1.8&2.4&4.2&8.2&12.9&10533&3&6&9&22&26\\
		&\hyperlink{target:Disjunctive}{$Disj(\epsilon^*)$}&4.3&2.0&2.5&4.4&12.9&16.2&10&4&6&10&25&29\\
		&\hyperlink{target:anytime}{$Anyt(\epsilon^*)$}&4.5&1.8&2.5&4.4&9.3&16.1&11&3&6&10&22&28\\
	\end{tabular*}      
\end{table}
\clearpage
\bibliographystyle{apalike}
\bibliography{sn-bibliography}
\end{document}